\documentclass[12pt]{compositio}
\usepackage{amssymb, amscd, amsmath, float, graphicx, longtable,color,
  enumerate
}
\usepackage{mathrsfs}
\usepackage{psfrag}
\usepackage[all]{xy}
\SelectTips{eu}{} 
\xyoption{curve}

\newcommand{\hide}[1]{}

\usepackage[colorlinks]{hyperref}

\DeclareMathOperator{\charac}{char}

\DeclareMathOperator{\coh}{{{coh}}}

\DeclareMathOperator{\cone}{{\mathsf{cone}}}

\DeclareMathOperator{\diag}{diag}
\DeclareMathOperator{\End}{End}
\DeclareMathOperator{\Ext}{Ext}

\DeclareMathOperator{\GL}{GL}

\DeclareMathOperator{\Grass}{Grass}

\DeclareMathOperator{\Hom}{Hom}
\DeclareMathOperator{\Id}{\operatorname{Id}}

\DeclareMathOperator{\ind}{ind}
\DeclareMathOperator{\im}{im}

\DeclareMathOperator{\rad}{rad}

\DeclareMathOperator{\Rep}{\operatorname{Rep}}

\DeclareMathOperator{\Sym}{Sym}

\renewcommand{\phi}{\varphi}

\renewcommand{\bar}{\overline}

\renewcommand{\mod}{\operatorname{mod}}
\renewcommand{\leq}{\leqslant}
\renewcommand{\le}{\leqslant}
\renewcommand{\geq}{\geqslant}
\renewcommand{\ge}{\geqslant}

\renewcommand{\to}{\longrightarrow}

\newcommand\into{{\hookrightarrow}}
\newcommand\xto{\ -\negmedspace\negmedspace\negthinspace\xrightarrow}
\newcommand\onto{\twoheadrightarrow}

\newcommand{\Wedge}{{\textstyle\bigwedge}}

\renewcommand{\L}{{\textstyle\bigwedge}}

\newcommand{\svee}{{\scriptscriptstyle \vee}}

\newcommand{\GG}{{\mathbb G}}

\newcommand{\PP}{{\mathbb P}}

\newcommand{\ZZ}{{\mathbb Z}}

\newcommand{\cala}{{\mathcal A}}
\newcommand{\calb}{{\mathcal B}}
\newcommand{\calc}{{\mathcal C}}
\newcommand{\cald}{{\mathcal D}}

\newcommand{\calf}{{\mathcal F}}

\newcommand{\calk}{{\mathcal K}}
\newcommand{\call}{{\mathcal L}}

\newcommand{\calo}{{\mathcal O}}

\newcommand{\calq}{{\mathcal Q}}
\newcommand{\calr}{{\mathcal R}}
\newcommand{\cals}{{\mathcal S}}
\newcommand{\calt}{{\mathcal T}}
\newcommand{\calu}{{\mathcal U}}

\newcommand{\calox}{{\calo_X}}

\newcommand{\gr}{\operatorname{gr}}
\newcommand{\rep}{\operatorname{rep}}
\def\uRHom{\operatorname {{\mathrm R}\mathcal{H}\!\mathit{om}}}

\def\cHom{{\mathcal H}\!om}
\newcommand{\rHom}[1][{}]{{\mathrm R}^{#1}\!\Hom}

\theoremstyle{plain}
\newtheorem{theorem}{Theorem}

\newtheorem{prop}[theorem]{Proposition}
\newtheorem{proposition}[theorem]{Proposition}

\newtheorem{lemma}[theorem]{Lemma}

\newtheorem{cor}[theorem]{Corollary}
\newtheorem{corollary}[theorem]{Corollary}

\newtheorem*{conjecture*}{Conjecture}
\newtheorem*{theorem*}{Theorem}
\newtheorem*{prop*}{Proposition}

\theoremstyle{definition}

\newtheorem{definition}[theorem]{Definition}

\newtheorem{remark}[theorem]{Remark}

\newtheorem{example}[theorem]{Example}

\numberwithin{theorem}{section}
\numberwithin{equation}{section}

\let\oldmarginpar\marginpar
\def\marginpar#1{\oldmarginpar{\tiny #1}}

\begin{document}

\title[The derived category of Grassmannians]{%
On the derived category of Grassmannians  in arbitrary characteristic
}

\author[R.-O. Buchweitz]{Ragnar-Olaf Buchweitz}
\email{ragnar@utsc.utoronto.ca}\address{Dept.\ of Computer and Mathematical Sciences, University of
Toronto Scarborough, Toronto, Ont.\ M1C 1A4, Canada}

\author[G.J. Leuschke]{Graham J. Leuschke}
\email{gjleusch@math.syr.edu}
\address{Dept.\ of Mathematics, Syracuse University,
Syracuse NY 13244, USA}

\author[M. Van den Bergh]{Michel Van den Bergh}
\email{michel.vandenbergh@uhasselt.be}
\address{Departement WNI, Universiteit Hasselt, 3590
  Diepenbeek, Belgium}

\thanks{The first author was partly supported by NSERC grant
  3-642-114-80. The second author was partly supported by NSF grant
  DMS~0902119.  The third author is director of research at the FWO\@.
  Part of this research was supported through the programme ``Research
  in Pairs'' by the Mathematisches Forschungsinstitut Oberwolfach in
  July 2012, and the work was completed at MSRI.  We thank both
  institutions for their hospitality.}

\date{\today}



\begin{abstract}
  In this paper we consider Grassmannians in arbitrary characteristic.
  Generalizing Kapranov's well-known characteristic-zero results we
  construct dual exceptional collections on them (which are however
  not strong) as well as a tilting bundle. We show that this tilting
  bundle has a quasi-hereditary endomorphism ring and we identify the
  standard, costandard, projective and simple modules of the latter.   
\end{abstract}
\maketitle

\section{Introduction}
Throughout $K$ is a field of arbitrary characteristic. Let $X$ be a
smooth algebraic variety over $K$ and let $\cald$ be its bounded
derived category of coherent sheaves.  An object $\calt\in \cald$ is
called a \emph{tilting object} if it classically generates $\cald$
(i.e.\ the smallest thick subcategory of $\cald$ containing $\calt$ is
$\cald$ itself) and $\Hom_{\calox}(\calt,\calt[i])=0$ for $i\neq 0$.

If $\calt$ is a tilting object in $\cald$ and $A=\End_{\calox}(\calt)$
then the functor $\rHom_{\calox}(\calt,-)$ defines an equivalence
$\cald\cong D^b(\mod A^\circ)$. If in addition $\calt$ is a vector
bundle then we call $\calt$ a \emph{tilting bundle}.

\medskip

A sequence of objects $E_1,E_2,\ldots,E_d$ which classically generates
$\cald$ is called an \emph{exceptional sequence} if
$\rHom_{\calox}(E_j,E_i)=0$ for $j>i$ and
$\rHom_{\calox}(E_i,E_i)=K$. An exceptional sequence is \emph{strongly
  exceptional} if in addition $\Ext^k_{\calox}(E_i,E_j)=0$ for all
$i,j$ and $k\neq 0$. Obviously if $(E_i)_i$ is strongly exceptional then
$\calt=\bigoplus_i E_i$ is a tilting object in $\cald$.

Two exceptional sequences $E_1,E_2,\ldots,E_d$ and
$F_d,F_{d-1},\ldots,F_1$ are said to be \emph{dual} if
\[
\rHom_{\calox}(E_i,F_j)=\delta_{i,j}\cdot K\,.
\]

We now specialize to the the case where $X$ is the Grassmannian $\GG =
\Grass(l,F) \cong \Grass(l,m)$ of $l$-dimensional subspaces of an
$m$-dimensional $K$-vector space $F$.  On $\GG$ we have a tautological
exact sequence of vector bundles
\begin{equation}
\label{ref-1.1-0}
0 \to \calr \to F^\svee \otimes_K \calo_\GG \to \calq \to 0
\end{equation}
in which $\calq$ has rank $l$ and $\calr$ has rank $m-l$.  When $K$ is
a field of characteristic zero, Kapranov~\cite{Kapranov:1988}
constructs a pair of dual strongly exceptional sequences on $\GG$
which we now describe.  For a partition~$\alpha$ let~$L^\alpha$ be the
associated Schur functor; our conventions are that $L^{(t)}V=\Sym^t V$
and $L^{(1^t)}V=\Wedge^t V$.  Further let $\alpha'$ be the transpose
of $\alpha$ and let $|\alpha|=\sum_i \alpha_i$ be its degree.
\begin{theorem}
[(Kapranov \cite{Kapranov:1988})] 
\label{thm:kapranov}
  Assume that $K$ has
  characteristic zero.  Let $B_{u,v}$ be the set of partitions with at
  most $u$ rows and at most~$v$ columns equipped with a total ordering
  $\prec$ such that if $|\alpha|<|\beta|$ then $\alpha\prec
  \beta$. Let $\bar{B}_{u,v}$ be the same as $B_{u,v}$ but equipped
  with the opposite ordering. Then there are strongly exceptional
  sequences on $\GG$ given by
  \[
  (L^\alpha \calq)_{\alpha\in B_{l,m-l}}
  \quad\text{and}\quad
  (L^{\alpha'} \calr)_{\alpha\in \bar{B}_{l,m-l}}\,.
  \]
  In particular the vector bundle
  \[
  \calk = \bigoplus_{\alpha \in B_{l,m-l}} L^\alpha \calq
  \]
  is a tilting bundle on $\GG$.  Moreover the exceptional sequences
  $(L^\alpha \calq)_{\alpha\in B_{l,m-l}}$ and $(L^{\alpha'}
  \calr[|\alpha|])_{\alpha\in \bar{B}_{l,m-l}}$ are dual.
\end{theorem}
For $K$ a field of positive characteristic $p$,
Kaneda~\cite{Kaneda:2008} shows that $\calk$ remains tilting as long
as $p \geq m-1$.  However $\calk$ fails to be tilting in very small
characteristics.  

\begin{example}
  \label{ref-4.1-17}
  Assume that $K$ has characteristic $2$ and put $\GG = \Grass(2,4)$.
  Then the short exact sequence
  \begin{equation}\label{ref-4.3-18}
  0 \to \L^2 \calq \to \calq \otimes \calq \to \Sym^2 \calq \to 0
  \end{equation}
  is non-split.  This follows for example from
  Theorem~\ref{ref-6.4-28} below and the fact that the sequence of
  $\GL(2)$-representations
  \[
  0\to \L^2 V\to V\otimes V\to \Sym^2  V\to 0
  \]
  is not split, where $V=K^2$ is the standard representation.  In
  particular $\Ext_{\calo_\GG}^1(\Sym^2\calq,\L^2\calq) \neq 0$, so that
  $\Sym^2\calq$ and $\L^2 \calq$ are not common direct summands of a
  tilting bundle on $\GG$ in characteristic two.
\end{example}

\medskip

In this note we give a tilting bundle on $\GG$ which exists in
arbitrary characteristic.  For a partition
$\alpha=[\alpha_1,\ldots,\alpha_p]$ put
\[
\L^\alpha \calq=\L^{\alpha_1}\calq\otimes_{\GG}\cdots\otimes_{\GG} 
\L^{\alpha_p}\calq\,.
\]
Our first main theorem is the following.
\begin{theorem}  \label{ref-1.2-1}
  Define a vector bundle on $\GG$ by
  \begin{equation}
    \label{ref-1.2-2}
    \calt = \bigoplus_{\alpha\in B_{l,m-l}} \L^{\alpha'} \calq\,.
  \end{equation}
  Then $\calt$ is a tilting bundle on $\GG$.
\end{theorem}
In characteristic zero we recover Kapranov's tilting bundle, up to
multiplicities, by working out the tensor products
in~\eqref{ref-1.2-2} using Pieri's formula.

\medskip

The proof of Theorem~\ref{ref-1.2-1} depends on the following
vanishing result which we will also use in \cite{BLVdB2}.
\begin{prop}
  \label{ref-1.3-3} 
  For $\alpha\in B_{l,m-l}$ and  $\beta$ an arbitrary partition we have for $i>0$
  \begin{equation}
    \label{ref-1.3-4}
    \Ext^i_{\calo_\GG}(\L^{\alpha'}\calq,L^\beta\calq)=0\,.
  \end{equation}
  Furthermore if $|\beta| < |\alpha|$ then we have as well
  $\Hom_{\calo_\GG}(\L^{\alpha'}\calq,L^\beta\calq)=0$. 
\end{prop}
In our next result we show that Kapranov's characteristic-zero result
can be partially salvaged in arbitrary characteristic. 
\begin{theorem}
  [(see Theorem~\ref{ref-8.4-47} below)]
 \label{ref-1.4-5}
  There exists a total ordering $\prec$ on $B_{l,m-l}$ such that
  \[
  (L^\alpha \calq)_{\alpha\in B_{l,m-l}}
\quad\text{and}\quad
  (L^{\alpha'} \calr[|\alpha|])_{\alpha\in \bar{B}_{l,m-l}}
  \]
  are dual exceptional collections on $\GG$, where $\bar{B}_{l,m-l}$ is
  $B_{l,m-l}$ equipped with the opposite ordering.
\end{theorem}
We use this result to obtain another proof of Kaneda's result that
$\calk$ remains tilting in characteristics $p \geq m-1$
(Corollary~\ref{cor:kaneda}).

The proof of Theorem~\ref{ref-1.4-5} goes through the construction of
a nice semi-orthogonal decomposition \cite{Bondal-Kapranov:1989} on $D^b(\coh(\GG))$
which we summarize in the following theorem.
\begin{theorem}
  [(see Theorem~\ref{ref-6.6-31} below)] 
  \label{thm:semiorth}
  There is a
  semi-orthogonal decomposition 
  \[
  D^b(\coh(\GG)) = \left\langle \cald_0,\ldots,\cald_{l(m-l)}\right\rangle
  \]
  where for $d = 0, \dots, l(m-l)$, $\cald_d$ is the derived category
  of the generalized Schur algebra associated to the representations
  whose composition factors have highest weight $\alpha\in B_{l,m-l}$
  such that $|\alpha|=d$.
\end{theorem}

The connection between Theorem~\ref{ref-1.2-1} and~\ref{ref-1.4-5}
depends on the theory of quasi-hereditary
algebras~\cite{Dlab-Ringel:1992}. In
this regard we have the following additional result.
\begin{theorem}\label{thm:qhed}
  Let $\calt$ be as in Theorem~\ref{ref-1.2-1} and put
  $A=\End_{\calo_{\GG}}(\calt)$.  Then $A$ is
  quasi-hereditary. Furthermore the homogeneous bundles
  $(L^\alpha\calq)_{\alpha\in B_{l,m-l}}$ correspond to the standard
  right $A$-modules and $(L^{\alpha'}\calr[|\alpha|])_{\alpha\in
    \bar{B}_{l,m-l}}$ correspond to the costandard right $A$-modules.
\end{theorem}
This theorem is a special case of Theorem~\ref{ref-9.1-49} below in
which we also characterize the simple and projective right
$A$-modules.

\medskip
The authors wish to thank Vincent Franjou, Catharina Stroppel and
Antoine Touz\'e for help with references.  After this paper was posted
to the arXiv, we learned that A.~Efimov has obtained
results similar to Theorems~\ref{thm:semiorth} and~\ref{thm:qhed} by
different methods.

\section{Some preliminaries on representation theory}
\label{sec:prelims}
Throughout we use \cite{Jantzen} as a convenient reference for facts
about algebraic groups.  If $H\subset G$ is an inclusion of algebraic
groups over the ground field $K$, then the restriction functor from
rational $G$-modules to rational $H$-modules has a right adjoint
denoted by $\ind^G_H$~(\cite[I.3.3]{Jantzen}).  Its right derived
functors are denoted by $R^i\ind^G_H$. For an inclusion of groups
$K\subset H\subset G$ and $M$ a rational $K$-representation there is a
spectral sequence \cite[I.4.5(c)]{Jantzen}
\begin{equation}
  \label{ref-3.1-7}
  E^{pq}_2 \colon
  R^p\ind_H^G  R^q\ind_{K}^H M \Longrightarrow R^{p+q}\ind_K^GM\,.
\end{equation}
If $H \subset G$ are group schemes such that $G/H$ is a scheme, and
$V$ is a finite-dimensional $H$-representation, then $\call_{G/H}(V)$
is by definition the $G$-equivariant vector bundle on $G/H$ given by
the sections of $(G\times V)/H$. The functor $\call_{G/H}(-)$ defines
an equivalence~\cite[Theorem 2.7]{Cline-Parshall-Scott:1983} between
the finite-dimensional $H$-representations and the $G$-equivariant
vector bundles on $G/H$. The inverse of this functor is given by
taking the fiber in $[H] \in G/H$.

If $G/H$ is a scheme then $R^n\ind_H^G$ may be computed as
\cite[Prop.\ I.5.12]{Jantzen}
\begin{equation}
  \label{ref-3.2-8}
  R^n\ind^G_H M=H^n(G/H,\call_{G/H}(M))\,.
\end{equation}

We now assume that $G$ is a split reductive group with a given split
maximal torus and Borel $T\subset B\subset G$. We let $X(T)$ be the
character group of $T$ and we identify the elements of $X(T)$ with the
one-dimensional representations of $T$.  The set of roots (the weights
of $G$ on $\operatorname{Lie}{G}$) is denoted by $R$. We have $R=R^-\coprod
R^+$ where the negative roots $R^-$ represent the roots of
$\operatorname{Lie} B$. For $\alpha\in R$ we denote the corresponding
coroot in $Y(T)=\Hom(X(T),\ZZ)$ \cite[II.1.3]{Jantzen} by
$\alpha^\svee$.  The natural pairing between $X(T)$ and $Y(T)$ is
denoted by $\left\langle -,-\right\rangle$.  A weight $\lambda\in X(T)$ is
dominant if $\left\langle \lambda,\alpha^\svee\right\rangle\ge 0$ for all positive
roots $\alpha$. The set of dominant weights is denoted by $X(T)_+$.
The set $X(T)$ is naturally partially ordered by putting $\lambda\le
\mu$ if $\mu-\lambda$ is a sum of positive roots.

The following is the celebrated Kempf vanishing result
(\cite{Kempf:1976}, see also~{\cite[II.4.5]{Jantzen}}).
\begin{theorem}
  \label{ref-3.1-9}  
  If $\lambda\in X(T)_+$ then
  $R^i\ind_B^G\lambda = H^i(G/B,\call_{G/B}(\lambda))$
  vanishes for all strictly positive $i$.\qed
\end{theorem}


We now restrict to $G=\GL(m)$. In this case we let $T$ be the diagonal
matrices in $G$ and $B$ the lower triangular matrices.  The weights of
$T$ can be identified with $m$-tuples of integers
$[\alpha_1,\ldots,\alpha_m]$ via $\diag(z_1,\ldots,z_m)\mapsto
z_1^{\alpha_1}\cdots z_m^{\alpha_m}$. Thus $X(T)\cong Y(T)\cong
\ZZ^m$. Under this identification roots and coroots coincide and are
given by $(0,\ldots,0,\pm 1,0,\ldots,0, \mp1,0,\ldots,0)$. The pairing
between $X(T)$ and $Y(T)$ is the standard Euclidean scalar product and
hence $X(T)_+=\{[\alpha_1,\ldots,\alpha_m]\mid \alpha_i\ge
\alpha_j\text{ for $i\le j$}\}$. A dominant weight with only
non-negative entries will be called a \emph{partition}. Mentally we
represent a partition by its Young diagram, with the length of the
rows corresponding to the entries.  The sum $\sum_i \alpha_i$ is the
\emph{degree} of the weight $\alpha$ and is denoted by $|\alpha|$. We
say that a representation has degree $d$ if all its weights have
degree $d$. We say that a representation is \emph{polynomial} if all
its weights contain only non-negative entries.

If $\alpha = [\alpha_1, \dots, \alpha_m]$ is a partition then we
denote by $L^\alpha$, $K^\alpha$ the corresponding \emph{Schur} and
\emph{Weyl functors}. More precisely for a vector space (or a vector
bundle) $V$ define for a partition $\alpha$
\[
\L^\alpha V=\bigotimes_i \L^{\alpha_i} V\qquad
\Sym^\alpha V=\bigotimes_i \Sym^{\alpha_i} V\qquad
D^\alpha V=\bigotimes_i D^{\alpha_i} V
\]
where in particular $D^u V=(V^{\otimes u})^{S_u}$ is the $u^\text{th}$
divided power representation, $\Wedge^u V$ is the exterior power, and
$\Sym^u V$ is the symmetric power.

Then we put with $d=|\alpha|$:
\begin{align}
  L^\alpha V&=\im\left(\L^{\alpha'} V\xto{a}
    V^{\otimes d}\xto{s}  \Sym^{\alpha} V\right)\label{ref-3.3-10}\\
  K^\alpha V&=\im\left(D^\alpha V\xto{s}
    V^{\otimes d}\xto{a}  \L^{\alpha'} V\right)\label{ref-3.4-11}
\end{align}
where $a$ and $s$ are respectively the anti-symmetrization map and the
symmetrization map. Their precise form is obtained from a filling of
the Young diagram associated to $\alpha$
(see~\cite[\S8.1]{Fulton:1997}). The resulting representations
$K^\alpha V$, $L^\alpha V$ are independent of this labeling.

In the sequel we freely pass between the functor point of view and the
representation theory point of view using the following lemma.  If
$\lambda\in X(T)_+$ then
$H^0(\lambda)\overset{\text{def}}{=}\ind_B^G\lambda$ is a so-called
\emph{induced representation} with highest weight $\lambda$.  Dually
one defines the corresponding \emph{Weyl representation} as
$V(\lambda)=H^0(-w_0\lambda)^\svee$ where $w_0$ is the longest element
of the Weyl group \cite[\S 2.13]{Jantzen}.
\begin{lemma} 
  Let $V$ be the standard representation of $\GL(m)$ and let $\alpha$
  be a partition.  Then
  \begin{align}
    L^\alpha V&=H^0(\alpha)\label{ref-3.5-12}\\
    K^\alpha V&=V(\alpha)\label{ref-3.6-13}
  \end{align}
\end{lemma}
\begin{proof}
  The identity \eqref{ref-3.5-12} is
  \cite[(4.1.10)]{Weyman:book}\footnote{Note that our $L^\alpha$ is
    $L_{\alpha'}$ in \cite{Weyman:book}.}. To prove
  \eqref{ref-3.6-13} we note that by \cite[II.2.13(2)]{Jantzen} we
  have $V(\alpha)={}^\tau\, H^0(\alpha)$, where ${}^\tau M$ is
  $M^\svee$ as a vector space and $g\in G$ acts on $\phi\in M^\svee$
  via $g\cdot \phi= \phi\circ g^t$ where $(-)^t$ denotes
  transposition. Clearly $M\mapsto {}^\tau M$ is a contravariant monoidal
  functor and furthermore one verifies
  \begin{align*}
    {}^\tau \Sym^u V&=D^u V\\
    {}^\tau \L^u V&=\L^u V\,.
  \end{align*}
  Applying ${}^\tau\,(-)$ to the right-hand side of \eqref{ref-3.3-10}
  yields the right-hand side of \eqref{ref-3.4-11}, finishing the
  proof.
\end{proof}
According to \cite[Prop.\ II.2.4]{Jantzen} $L^\alpha V$ has a simple
socle which we denote by $\Sigma^\alpha$. According to \cite[\S
II.2.6]{Jantzen} $K^\alpha V$ has a simple top, which is also equal to
$\Sigma^\alpha$.  (Recall that the socle is the sum of all simple
submodules, while the top of a module $M$ is $M/\rad M$.)

We also state for easy reference the following characteristic-free
versions of the Cauchy formula and the Littlewood-Richardson rule.
See~\cite[(2.3.2), (2.3.4)]{Weyman:book}. 

\begin{theorem}[(Boffi~\cite{Boffi:1988},
  Doubilet-Rota-Stein~\cite{Doubilet-Rota-Stein:1974})]
  \label{thm:cauchyLR}
  Let $V$ and $W$ be $K$-vector spaces and let $\alpha$ and $\beta$ be
  dominant weights.
  \begin{enumerate}[\quad(i)]
  \item \label{item:cauchy} There is a natural filtration on $\Sym^t(V
    \otimes W)$ whose associated graded object is a direct sum with
    summands tensor products $L^\gamma V \otimes L^{\gamma'} W$ of Schur
    functors.
  \item \label{item:LR} There is a natural filtration on $L^\alpha V
    \otimes L^\beta V$ whose associated graded object is a direct sum
    of Schur functors $L^\gamma V$. The $\gamma$ that appear, and
    their multiplicities, can be computed using the usual
    Littlewood-Richardson rule.
  \end{enumerate}
  If $\charac K=0$ then the filtrations above degenerate to direct
  sums. 
\end{theorem}

\section{Proofs of Theorem~\ref{ref-1.2-1} and Proposition~\ref{ref-1.3-3}}
We stick to the notation already introduced in the introduction.  We
will identify $\GG=\Grass(l,F)$ with $\Grass(m-l,F^\svee)$ via the
correspondence $(V\subset F)\mapsto ((F/V)^\svee\subset F^\svee)$.

For convenience we choose a basis  $(f_i)_{i=1,\ldots,m}$ for $F$
and a corresponding dual basis $(f^\ast_i)_{i}$ for $F^\svee$. We view
$\GG$ as the homogeneous space $G/P$ with $G = \GL(F^\svee)=\GL(m)$ and
$P\subset G$ the parabolic subgroup stabilizing the point $(W\subset
F^\svee)\in \GG$ where $W=\sum_{i=l+1}^{m} K f_i^\ast$.  As above let
$T$ and $B$ be respectively the diagonal matrices and the lower
triangular matrices in $G$.
\[
\includegraphics[scale=0.75]{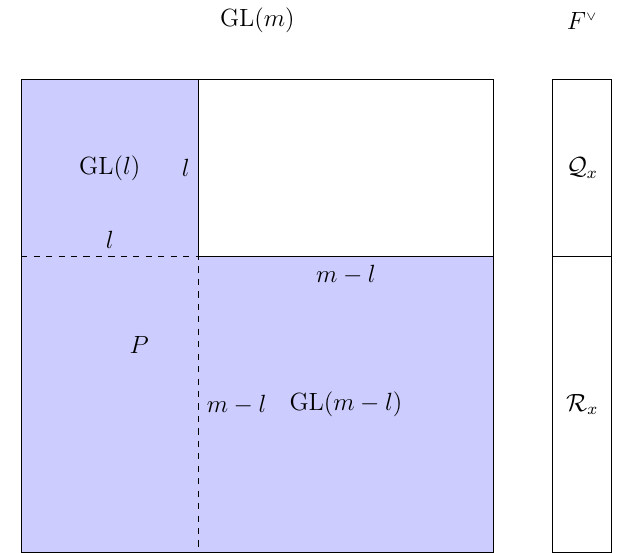}
\]

Let $H=G_1\times G_2=\GL(l)\times \GL(m-l)\subset \GL(m)$ be the
Levi-subgroup of $P$ containing~$T$. We put $B_i=B\cap G_i$ and
$T_i=T\cap G_i$.  We denote the standard representations of $G_1$ and
$G_2$ by $V$ and $W$ respectively. Thus for $x=[P]\in G/P$ we have
$V=\calq_x$ and $W=\calr_x$; equivalently $\calq = \call_\GG(V)$ and
$\calr = \call_\GG(W)$.  (Throughout we silently view
$G_i$-representations as $P$-representations to apply $\call_\GG(-)$.)
It follows that $\Wedge^{\alpha'}\calq = \call_\GG(\Wedge^{\alpha'}V)$
and $L^\alpha \calq = \call_\GG(L^\alpha V)$.

For use in the proof below we fix an additional parabolic $P^\circ$ in
$G$ given by the stabilizer of the flag $(\sum_{i\ge p} K
f^\ast_i)_{p=1,\ldots,l}$.  We let $G^\circ=\GL(m-l+1)\subset
P^\circ\subset G=\GL(m)$ be the lower right $(m-l+1\times
m-l+1)$-block in $\GL(m)$. We put $T^\circ=T\cap G^\circ$ and
$B^\circ=B\cap G^\circ$, i.e.~$B^\circ$ is the set of lower
triangular matrices in $G^\circ$ and $T^\circ$ is the set of diagonal
matrices.
\begin{figure*}[htpb]
\[
\includegraphics[scale=0.75]{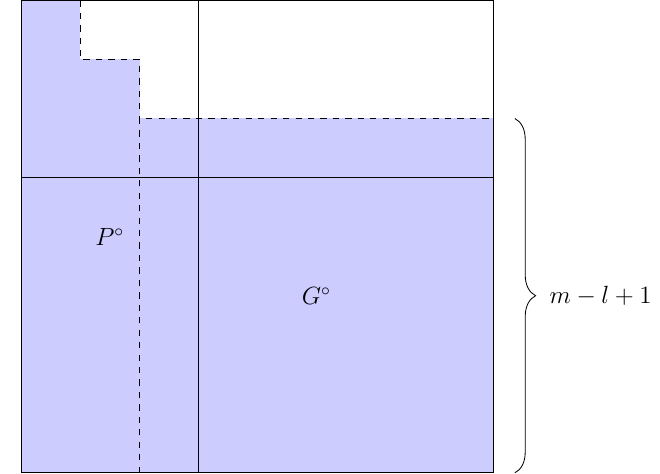}
\]
\end{figure*}

\begin{proof}[Proof of Proposition~\ref{ref-1.3-3}]
  The usual spectral sequence argument implies that
  $\Ext^i_{\calo_\GG}(\L^{\alpha'}\calq, L^\beta\calq)$ is the
  $i^\text{th}$ cohomology of $\cHom_{\calo_\GG}(\L^{\alpha'}\calq,
  L^\beta\calq) \cong (\L^{\alpha'}\calq)^\svee \otimes_\GG L^\beta\calq$,
  so we must show
  \[
  H^i(\GG,\L^{u_1}\calq^\svee\otimes_\GG \cdots
  \otimes_\GG\L^{u_{m-l}}\calq^\svee\otimes_\GG L^\beta\calq)=0
  \]
  for all $i>0$ and $u_1\ge\cdots\ge u_{m-l}\ge 0$, and also for $i=0$
  if $\sum u_i > |\beta|$.

  Using the identity
  \[
  \left(\L^u\calq\right)^\svee = \L^{l-u} \calq \otimes \left(\L^l
    \calq\right)^\svee
  \]
  and the fact (Theorem~\ref{thm:cauchyLR}(\ref{item:LR})) that
  $\Wedge^{l-u}\calq \otimes L^\beta \calq$ is filtered by
  subquotients of the form $L^\gamma \calq$ with $|\gamma| \geq
  |\beta|$, we can reduce immediately to the case $u_1 = \cdots =
  u_{m-l} = l$, replacing $\beta$ with a larger partition if
  necessary.  The tautological exact sequence~\eqref{ref-1.1-0} allows
  us to write
  \[
  \left(\L^l \calq\right)^\svee = \L^m F \otimes_K
  \L^{m-l}\calr\,.
  \]
  Thus it is enough to prove that for an arbitrary partition $\gamma$,
  \[
  L^\gamma\calq \otimes_\GG \L^{m-l}\calr \otimes_\GG \cdots \otimes_\GG
  \L^{m-l}\calr
  \]
  (with $m-l$ factors of $\L^{m-l}\calr$) has vanishing higher
  cohomology.  Using \eqref{ref-3.2-8} we see that we must prove that
  for $i >0$ we have
  \begin{equation}
    \label{ref-4.1-15}
    R^i \ind_P^{G} \left(L^\gamma V \otimes \L^{m-l}W
    \otimes \cdots \otimes \L^{m-l}W\right) =0\,, 
  \end{equation}
  where as above $V$, $W$ are the standard representations of $G_1$,
  $G_2$.  Since $V$ has rank $l$, we may assume that $\gamma$ has
  at most~$l$ entries.  Put
  $\bar{\gamma}=[\gamma_1,\ldots,\gamma_l,m-l,\ldots,m-l]\in X(T)$.
  We have
  \[
  \begin{aligned}
    L^\gamma V \otimes \L^{m-l}W \otimes \cdots
    \otimes \L^{m-l}W 
    &= \ind_{B_1}^{G_1} L^\gamma V \otimes \ind_{B_2}^{G_2}
    L^{(m-l)^{m-l}}W \\
    &= \ind_B^P \bar{\gamma}\,.
  \end{aligned}
  \]
  It is clear that $\bar\gamma$ is dominant when viewed as a weight
  for $T$ considered as a maximal torus in $H=G_1\times G_2$ with
  respect to the Borel subgroup $B_1\times B_2$.  So Kempf vanishing
  implies that $R^i \ind_B^P\bar\gamma=R^i\ind_{B_1\times
    B_2}^{G_1\times G_2}\bar\gamma =0$ for all $i >0$.

  Thus the spectral sequence~\eqref{ref-3.1-7} degenerates and we
  obtain
  \begin{equation}
    \label{ref-4.2-16}
    R^i \ind_P^G \left(L^\gamma V\otimes \L^{m-l}W
      \otimes \cdots \otimes \L^{m-l}W\right)
    =
    R^i \ind_B^{G} \bar{\gamma}\,.
  \end{equation}
  Thus if $\bar \gamma$ is dominant (i.e.\ $\gamma_l \geq m-l$) then
  the desired vanishing~\eqref{ref-4.1-15} follows by invoking Kempf
  vanishing again.

  Assume then that $\bar \gamma$ is not dominant, i.e.\ $0 \leq
  \gamma_l < m-l$. We claim that $R^i\ind^{P^\circ}_B \bar{\gamma}=0$
  for all $i$. Then by the spectral sequence~\eqref{ref-3.1-7} applied
  to $B\subset P^\circ \subset G$ we obtain that $R^i\ind^{G}_B
  \bar{\gamma}=0$ for all $i$.

  To prove the claim we note that $P^\circ/B\cong G^\circ/B^\circ$ and
  hence by \eqref{ref-3.2-8} $R^i\ind^{P^\circ}_B
  \bar{\gamma}=R^i\ind^{G^\circ}_{B^\circ} (\bar{\gamma}\mid
  T^\circ)$. In other words we have reduced ourselves to the case
  $l=1$ (replacing $m$ by $m-l+1$).

  So now we assume $l=1$. Thus $\GG=\PP(F) \cong \PP^{m-1}$, which we
  write as $\PP$ for short.  The
  partition $\gamma$ consists of a single entry $\gamma_1$ and we have
  $\bar{\gamma}=[\gamma_1,m-1,\ldots,m-1]$.  Under the assumption
  $\gamma_1<m-1$ we have to prove $R^i\ind^{G}_B \bar{\gamma}=0$ for
  all $i$. Applying \eqref{ref-4.2-16} in reverse this means we have
  to prove that
  \[
  \calq^{\otimes \gamma_1}\otimes_\PP  \left(\L^{m-1}\calr\right)^{\otimes m-1}
  \]
  has vanishing cohomology on $\PP$. 

  We now observe $\calq\cong\calo_{\PP}(1)$ and since
  \[
  \calr \cong \ker (\calo_{\PP}^m\to \calo_{\PP}(1))
  \]
  we also have
  \[
  \L^{m-1}\calr \cong \calo_{\PP}(-1)
  \]
  so that
  \[
  \calq^{\otimes \gamma_1}\otimes_\PP \L^{m-1}\calr^{\otimes
    m-1} \cong \calo_{\PP}(-m+1+\gamma_1)\,.
  \]
  It is standard that this line bundle has vanishing cohomology when
  $0 \leq \gamma_1<m-1$, so we are done.

  For the last statement of the Proposition, observe that in the above
  argument if $|\beta|< |\alpha|$ then we are always in the case where
  $\bar\gamma$ is not dominant, and thus the vanishing holds also when
  $i=0$.
\end{proof}

\begin{proof}[Proof of Theorem~\ref{ref-1.2-1}]
  The main thing to prove is that $\Ext^i_{\calo_\GG}(\calt,\calt)=0$
  for $i \neq 0$. 
  Applying the characteristic-free
  Littlewood-Richardson rule Theorem~\ref{thm:cauchyLR}(\ref{item:LR}), we see that it suffices to prove
  that $ \calt^\svee\otimes_\GG L^\gamma\calq$ has vanishing higher
  cohomology whenever $\gamma$ is a partition with at most $l$ rows.
  This is the content of Proposition~\ref{ref-1.3-3}.

  Kapranov's resolution of the diagonal argument together with the
  characteristic-free version of Cauchy's
  formula (Theorem~\ref{thm:cauchyLR}(\ref{item:cauchy}))
  still implies that the vector bundle $\calk$ in
  Theorem~\ref{thm:kapranov} classically generates $D^b(\coh(\GG))$.
  See, for example,~\cite{Levine-Srinivas-Weyman:1989}.  Thus it
  suffices to show that $L^\alpha\calq$ for $\alpha\in B_{l,m-l}$ is
  in the thick subcategory $\calc$ generated by $\calt$.  Inductively,
  we may assume that $\alpha$ is such that $L^\beta \calq$ is in $\calc$
  for all $\beta$ less than $\alpha$ in the lexicographic ordering on
  partitions.

  Consider $\calu=\L^{\alpha'_1}\calq \otimes_\GG \cdots \otimes_\GG
  \L^{\alpha'_{l}}\calq$.  Then Pieri's formula, which is a special
  case of the Littlewood-Richardson rule, yields a filtration of
  $\calu$ with successive quotients $L^\beta \calq$ such that $\beta
  \leq \alpha$ and such that $L^\alpha \calq$ appears with
  multiplicity one.  Furthermore $\calu$ has a \emph{good filtration}
  \cite[\S II.4.16]{Jantzen}, one in which the $L^\beta \calq$
  appearing as quotients are in decreasing order for the lexicographic
  ordering on partitions, that is, the largest $\beta$ appear on
  top~\cite[II.4.16, Remark (4)]{Jantzen}.  Hence $\calu$ maps
  surjectively to $L^\alpha\calq$ and the kernel is an extension of
  various $L^\beta\calq$ with $\beta$ strictly smaller than $\alpha$
  in the lexicographic ordering. By the hypothesis all such
  $L^\beta\calq$ are in $\calc$. Since $\calu$ is in $\calc$ as well
  we obtain that $L^\alpha\calq$ is in $\calc$.
\end{proof}

\begin{remark}
  \label{ref-4.2-19}
  By \cite[Lemma (3.4)]{Donkin:1993} the indecomposable summands of
  $\Wedge^{\alpha'}V$ are precisely the \emph{tilting representations}
  for $\GL(V)$, so we can obtain the 
  following more economical tilting bundle for $\GG$
  \[
  \calt^\circ=\bigoplus_{\alpha\in B_{l,m-l}}\call_{\GG} (M^\alpha)\,,
  \]
  where $M^\alpha$ is the tilting $\GL(l)$-representation with highest
  weight $\alpha$~\cite[E.3]{Jantzen}.  Note however that the
  character of $M^\alpha$
  strongly depends on the characteristic.  Hence so does the nature of
  $\calt^\circ$.
\end{remark}

For use below we need  the following complement to
Proposition~\ref{ref-1.3-3}.
\begin{proposition}\label{ref-4.3}
  For every partition $\alpha$ and every polynomial
  $G_1$-representation $U$ of degree $<|\alpha|$ we have
  \begin{equation}
    \label{ref-4.4-20}
    \rHom_{\calo_{\GG}}(\L^{\alpha'}\calq,\call_{\GG}(U))=0\,.
  \end{equation}
\end{proposition}

\begin{proof} 
  It suffices to prove the claimed vanishing for $U$ simple of degree
  less than $|\alpha|$, so for $U=\Sigma^\beta$ with $\beta$ a partition
  such that $|\beta|<|\alpha|$. We do this by induction on
  $\beta$. Since $\Sigma^\beta$ is the socle of $L^\beta\calq$ we have
  by \cite[Prop.\ 6.15]{Jantzen} a short exact sequence
  \[
  0\to \Sigma^\beta\to L^\beta V \to S\to 0
  \]
  where $S$ is obtained through extensions involving only
  $\Sigma^\gamma$ with $\gamma<\beta$. By induction we may assume
  $\rHom_{\calo_{\GG}}(\L^{\alpha'}\calq,\call_{\GG}(S))=0$.  Then
  \eqref{ref-4.4-20} for $U=\Sigma^\beta$ follows from the final
  statement of Proposition~\ref{ref-1.3-3}.
\end{proof}

\section{Reminder on semi-orthogonal decompositions}\label{ref-5-22}
We recapitulate some facts concerning semi-orthogonal decompositions
that we need later. No originality is intended.

If $\cals$ is a triangulated category and $S$ is a set of objects then
we denote by $\left\langle S\right\rangle$ the smallest triangulated subcategory
of $\cals$ that contains $S$ and is closed under
isomorphisms. If $\cals=\left\langle S \right\rangle$ then we say that $S$
generates $\cals$ as a triangulated category.  (This is stronger than
``classically'' generating $\cals$ as in the Introduction.) 
\begin{definition}
  \label{ref-5.1-23}
  A \emph{semi-orthogonal decomposition} of a triangulated category
  $\cals$ is a sequence of full subcategories
  $\cala_1,\ldots,\cala_n\subset \cals$ generating $\cals$ as a
  triangulated category and such that $\Hom_\cals(\cala_j,\cala_i)=0$ for
  $j>i$. We denote such a semi-orthogonal decomposition by
  $\left\langle \cala_1,\ldots,\cala_n\right\rangle$. Sometimes we write
  $\cals=\left\langle \cala_1,\ldots,\cala_n\right\rangle$.
\end{definition}

If $X$ is an object in a triangulated category then a filtration $F$
of length $n$ on $X$ is a sequence of maps
\[
0=F_n X\to F_{n-1} X\to \cdots\to F_0X=X \,.
\]
We write $(\gr_F X)_i=\cone (F_{i+1}X\to F_i X)$.  The
following well-known lemma shows that Definition~\ref{ref-5.1-23} is
equivalent to the seemingly stronger one in \cite[Def. 2.3]{Kuznetsov:2009}.

\begin{lemma}
  Let $\left\langle\cala_1,\ldots,\cala_n\right\rangle$ be a semi-orthogonal decomposition of
  $\cals$.  Then every object $X$ in $\cals$ has a filtration $F$ of
  length $n$ such that $(\gr_F X)_i\in \cala_{i+1}$.
\end{lemma}

\begin{proof} By induction it is sufficient to prove this for
  $n=2$. In that case the result is \cite[Lemma 3.1]{Bondal:1989}.
\end{proof}

In order to work conveniently with semi-orthogonal decompositions one
needs a property called ``\emph{admissibility}''
\cite{Bondal-Kapranov:1989}.  
If $\cala$ is a  full triangulated
subcategory of a triangulated category $\cals$ then $\cala$ is (left,
right) admissible if the inclusion functor $\cala\to \cals$ has a
(left, right) adjoint, or equivalently if there exist semi-orthogonal
decompositions $\left\langle\cala,\cala'\right\rangle$ resp.\
$\left\langle\cala'',\cala\right\rangle$.  If $\cala$ is both left and
right admissible then it is said to be admissible.

A \emph{saturated} triangulated category is a $K$-linear triangulated
category $\cala$ such that for all $A,\ B\in \cala$ we have $\sum_i
\dim\Hom^i_\cala(A,B)<\infty$ and such that every co- or contravariant
cohomological functor $H^i\colon \cala\to \mod(K)$ satisfying
$\sum_i \dim H^i(A)<\infty$ is representable.  The derived category
of coherent sheaves on a smooth proper algebraic variety is a
particular example of a saturated triangulated category
\cite{BVdB,Bondal-Kapranov:1989}.

If $\cala$ is a saturated full triangulated subcategory of a $K$-linear
triangulated category $\cals$ then $\cala$ is admissible \cite[Prop
2.6]{Bondal-Kapranov:1989}.  Furthermore if $\cals$ is a saturated
triangulated category then every left/right admissible subcategory is
automatically admissible (and hence saturated).  This follows by
combining \cite[Prop.\ 2.6]{Bondal-Kapranov:1989} and \cite[Prop.\
2.8]{Bondal-Kapranov:1989}. From this we deduce that if we have a
semi-orthogonal decomposition
$\cals=\left\langle\cala_1,\ldots,\cala_n\right\rangle$ of a saturated
$\cals$ then all the ``slices''
$\left\langle\cala_i,\ldots,\cala_j\right\rangle$ are admissible and
saturated.

In particular if we put $\cals_{\le
  i}=\left\langle\cala_1,\ldots,\cala_i\right\rangle$ then this yields
a filtration $\cals_{\le 1}\subset \cdots \subset \cals_{\le n}=\cals$
by admissible subcategories. Let $\calb_i$ be the right orthogonal of
$\cals_{\le i-1}$ in $\cals_{\le i}$. Then we have semi-orthogonal
decompositions $\cals_{\le i} = \left\langle\calb_i,\cals_{\le
    i-1}\right\rangle$. Iterating we obtain a semi-orthogonal
decomposition
\[
\cals=\left\langle\calb_n,\calb_{n-1},\ldots,\calb_1\right\rangle
\]
such that
\[
\left\langle\cala_1,\ldots,\cala_i\right\rangle = \left\langle\calb_i,\ldots,\calb_1\right\rangle\,.
\]
This is called the semi-orthogonal decomposition (right) dual to
$\cals=\left\langle\cala_1,\ldots,\cala_n\right\rangle$. Note that \cite[(4)]{Kuznetsov:2009}
\begin{align*}
  \calb_i&=\cals_{\le i}\cap \cals^{\perp}_{\le i-1}\\
  &=\left\langle\cala_{i+1},\ldots,\cala_n\right\rangle^\perp\cap \left\langle\cala_{1},\ldots,\cala_{i-1}\right\rangle^\perp\\
  &=\left\langle\cala_1,\ldots,\cala_{i-1},\cala_{i+1},\ldots,\cala_n\right\rangle^\perp\,.
\end{align*}
In particular for $i\neq j$
\begin{equation}
  \label{ref-5.1-24}
  \Hom(\cala_i,\calb_j)=0\,.
\end{equation}
The following is also well-known \cite{Kuznetsov:2009}.
\begin{lemma}
  \label{ref-5.3-25} Assume $\cals$ is saturated. Let $\gamma_i$ be
  the composition of the canonical functors
  \[
  \gamma_i\colon \cala_i\to \cals_{\le i} \to \cals_{\le i}\big/\cals_{\le i-1}=\calb_i\,.
  \]
  Then $\gamma_i$ is an equivalence of categories. Furthermore we have for
  $A\in \cala_i$, $B\in \calb_i$
  \[
  \Hom_\cals(A,B)=\Hom_{\calb_i}(\gamma_i(A),B)= \Hom_{\cals}(\gamma_i(A),B)\,.
  \]
\end{lemma}

\begin{proof} We have semi-orthogonal decompositions
  \[
  \cals_{\le i}=\left\langle\cals_{\le i-1},\cala_i\right\rangle
  = \left\langle\calb_i,\cals_{\le i-1}\right\rangle\,.
  \]
  The fact that $\gamma_i$ is an equivalence follows from \cite[Lemma
  1.9]{Bondal-Kapranov:1989}.

  Let $\jmath\colon \cala_i\to \cals_{\le i}$ and  $\imath\colon \calb_i\to
  \cals_{\le i}$ be the inclusion functors and let $\imath^*$ be the left
  adjoint to $\imath$. Then $\gamma_i=\imath^\ast\circ \jmath$.  We have
  \begin{align*}
    \Hom_\cals(A,B)&=\Hom_{\cals_{\le i}}(\jmath A,\imath B)\\
    &=\Hom_{\cals_{\le i}}(\imath^\ast \jmath A,B)\\
    &=\Hom_{\calb_i}(\gamma_i(A),B)\,.
  \end{align*}
  The equality $\Hom_{\calb_i}(\gamma_i(A),B) = \Hom_{\cals}(\gamma_i(A),B)$ is
  just that $\calb_i$ is a full subcategory of $\cals$.
  \end{proof}

\section{Semi-orthogonal decompositions for Grassmannians}
In this section we write $\cald$ for the bounded derived category of
coherent sheaves on $\GG$. This is in particular a saturated category
(see \S\ref{ref-5-22}).  We will construct a semi-orthogonal
decomposition of $\cald$.

We start by observing that the proof of Theorem~\ref{ref-1.2-1}
actually shows
\begin{lemma}
  \label{ref-6.1-26}
  $\cald$ is generated as a triangulated category by $(\L^{\alpha'}
  \calq)_{\alpha\in B_{l,m-l}}$ (instead of just classically
  generated, see \S\ref{ref-5-22}). \qed
\end{lemma}

A set $S$ of dominant weights is \emph{saturated} if whenever
$\alpha\in S$ and $\beta<\alpha$ is dominant we have $\beta\in
S$. (Here and below ``$<$'' is the standard ordering on weights; see
\S\ref{sec:prelims}.)  The set $B_{l,m-l}$ is an example of a
saturated set for $\GL(l)$.  For $d\ge 0$ let $\calc_{d}$ be the category of
finite-dimensional $G_1=\GL(l)$-representations whose composition
factors have highest weights $\alpha$ satisfying $|\alpha|=d$ and
$\alpha\in B_{l,m-l}$.  Thus $\calc_d$ is a \emph{truncated category}
in the sense of \cite[Ch. A]{Jantzen} associated to a saturated set of
dominant weights. In particular $\calc_d$ is the category of finite modules
over a certain finite-dimensional algebra, called a \emph{generalized
  Schur algebra} \cite[\S A.16]{Jantzen}.

We collect some elementary facts about the derived category of
$\calc_d$.
\begin{lemma} Let $\Rep(G_1)$ be the category of rational
  $G_1$-representations and for each $d$ let
  $D^b_{\calc_d}(\Rep(G_1))$ be the bounded derived category of
  complexes of representations having cohomology in $\calc_d$.  The
  canonical functor
  \[
  D^b(\calc_d)\to D^b_{\calc_d}(\Rep(G_1))
  \]
  is an equivalence of categories. 
\end{lemma}

\begin{proof} That the functor is fully faithful follows from the fact
  that the Yoneda $\Ext$'s in $\calc_d$ are the same as those in the
  ambient category $\Rep(G_1)$ (see \cite[Prop.\ A.10]{Jantzen}).
  Essential surjectivity follows from fully faithfulness and the fact
  that the essential image contains the generating
  subcategory $\calc_d$.
\end{proof}

In the sequel we will simply confuse $D^b(\calc_d)$ and $D^b_{\calc_d}(\Rep(G_1))$.
\begin{lemma}
  \label{ref-6.3-27}
  The triangulated category $D^b(\calc_d)$ is generated by the
  representations $\Wedge^{\alpha'} V$ for $\alpha\in B_{l,m-l}$,
  $|\alpha|=d$, where as usual $V$ is the standard representation of
  $G_1$.
\end{lemma}

\begin{proof}
  This is of course well-known but for the convenience of the reader
  we give the proof.  Let $\cala$ be the full subcategory of
  $D^b(\calc_d)$ generated by $(\Wedge^{\alpha'} V)_{\alpha\in
    B_{l,m-l}, |\alpha|=d}$. It is sufficient to prove that $\cala$
  contains the simples $\Sigma^\alpha$ for $\alpha\in B_{l,m-l}$,
  $|\alpha|=d$.

  By reasoning similar to the proof of Theorem~\ref{ref-1.2-1} we
  see that $\cala$ contains $K^\alpha V$ for $\alpha\in B_{l,m-l}$
  with $|\alpha|=d$. By
  \cite[II.2.13]{Jantzen} $K^\alpha V$ has simple top $\Sigma^\alpha$
  and by the dual version of \cite[II.6.13]{Jantzen} the other
  Jordan-H\"older quotients of $K^\alpha V$ are of the form
  $\Sigma^\gamma$ with $|\gamma|=|\alpha|=d$ and $\gamma<\alpha$.
  Thus $\Sigma^\gamma\in \calc_d$. By induction we may assume that
  such $\Sigma^\gamma\in \cala$. Hence it follows that
  $\Sigma^\alpha\in \cala$.
\end{proof}

We define a functor
  \[
  \Phi_d \colon D^b(\calc_d)\to \cald
  \]
  by $\Phi_d( U) = \call_{\GG}(U)$ for $U \in \cald^b(\calc_d)$, where we
  view $U$ as a complex of  $P$-representations in the obvious way.

\begin{theorem}
  \label{ref-6.4-28}
  The functor
  $\Phi_d$
  is fully faithful.
\end{theorem}

\begin{proof} By Lemma~\ref{ref-6.3-27} it is sufficient to prove that
  for $\alpha,\beta\in B_{l,m-l}$ with  $|\alpha|=|\beta|=d$ the canonical
  map
  \begin{equation}
    \label{ref-6.1-29}
    \rHom_{G_1}(\Wedge^{\alpha'}V,\Wedge^{\beta'}V)\to 
    \rHom_{\calo_\GG}(\Wedge^{\alpha'}\calq,\Wedge^{\beta'}\calq)
  \end{equation}
  is an isomorphism (where we have used that $\call_{\GG}(V)=\calq$). Now
  $\Wedge^{\alpha'}V$ and $\Wedge^{\beta'}V$ are tilting
  representations \cite[Lemma (3.4)]{Donkin:1993} and so on the
  left-hand side of \eqref{ref-6.1-29} there are no higher
  $\Ext$'s. Likewise on the right-hand side there are no higher
  $\Ext$'s because of Proposition~\ref{ref-1.3-3}.

  So we only have to show that the map
  \[
  \Hom_{G_1}(\Wedge^{\alpha'}V,\Wedge^{\beta'}V)\to 
  \Hom_{\calo_\GG}(\Wedge^{\alpha'}\calq,\Wedge^{\beta'}\calq)
  \]
  is an isomorphism.  It is certainly injective: take the fiber of the
  right-hand side in $[P] \in G/P$ to get
  $\Hom_K(\Wedge^{\alpha'}V,\Wedge^{\beta'}V)$, and the composition is
  injective.  Thus we need only compute the $K$-dimensions of each
  side.  Since the Euler characteristic is independent of base
  field\footnote{This follows from the standard fact (see for
    example~\cite[3.8]{Danilov:1996}) that Euler characteristics are
    constant in families. To deduce characteristic independence one
    must use that everything in our setup can be defined over a
    commutative ring, and in particular over a discrete valuation ring
    of unequal characteristic.}%
  and the higher $\Ext$'s vanish, it suffices to do this in
  characteristic zero.

  Thus we assume that $K$ has characteristic zero. Then we may decompose
  $\Wedge^{\alpha'}V$, $\Wedge^{\beta'}V$ as sums of simple modules
  $L^\gamma V$, $L^\delta V$. Thus it is sufficient to prove that
  \[
  \Hom_{\GL(l)}(L^\alpha V,L^\beta V)\to \Hom_{\calo_\GG}(L^\alpha
  \calq,L^\beta\calq)
  \]
  is an isomorphism. Or in other words, since in characteristic zero
  the $L^\alpha V$ are simple,
  \[
  \Hom_{\calo_\GG}(L^\alpha \calq,L^\beta\calq)=\delta_{\alpha,\beta}\cdot K\,.
  \]
  This follows from the Littlewood-Richardson rule, combined with
  Bott's theorem (see e.g.\ \cite[\S 3.2, \S 3.3]{Kapranov:1988}).
\end{proof}

Now let $\cald_d$ be the essential image of $D^b(\calc_d)$ under
$\Phi_d$, i.e.\ the closure of that image under isomorphisms. From
Lemma~\ref{ref-6.3-27} we obtain:

\begin{corollary} 
  \label{ref-6.5-30}
  $\cald_d$ is generated by $\Wedge^{\alpha'}\calq$
  for $\alpha\in B_{l,m-l}$, $|\alpha|=d$. \qed
\end{corollary}

We have:
\begin{theorem}
  \label{ref-6.6-31}
  The triangulated category $\cald$ has a semi-orthogonal
  decomposition
  \begin{equation}
    \label{ref-6.2-32}
    \cald=  \left\langle \cald_0,\ldots,\cald_{l(m-l)}\right\rangle\,.
  \end{equation}
  Furthermore $\cald_{d}$ is generated by $L^{\alpha}\calq$ for
  $\alpha \in B_{l,m-l}$ with $|\alpha| = d$.
\end{theorem}

\begin{proof}  
  By Lemma~\ref{ref-6.1-26}, $\cald$ is generated by
  $\Wedge^{\alpha'}\calq$ for $\alpha\in B_{l,m-l}$. It follows from
  Corollary~\ref{ref-6.5-30} that $\cald$ is generated by $(\cald_d)_d$
  and that $\cald_{d}$ is generated by those $\L^{\alpha'}\calq$
  with $|\alpha| = d$.

  To complete the proof that \eqref{ref-6.2-32} is a semi-orthogonal
  decomposition we need that $\Hom(\cald_d,\cald_e)=0$ for $d>e$, or
  equivalently that
  $\rHom_\cald(\Wedge^{\alpha'}\calq,\Wedge^{\beta'}\calq)=0$ for
  $|\alpha|=d$, $|\beta|=e$.  This follows from
  Proposition~\ref{ref-4.3}. 
\end{proof}


The theorems we have stated have dual versions where $\calq$ is
replaced by $\calr$ and $B_{l,m-l}$ by $B_{m-l,l}$. We prove these by
passing to the dual Grassmannian $\Grass(m-l,F^\svee)$.

\begin{lemma} The vector bundle
  \begin{equation}
    \label{ref-6.3-33}
    \calt' = \bigoplus_{\alpha\in B_{m-l,l}} \Wedge^{\alpha'}\calr
  \end{equation}
  is a tilting bundle on $\GG$.
\end{lemma}

\begin{proof}
  Using the duality
  $\uRHom_{\calo_\GG}(-,\calo_{\GG})$ on $\cald$ it is sufficient to show
  that $\calt^{\prime \svee}$ is a tilting bundle. Now
  $\calt^{\prime\svee}$ is equal to $\bigoplus_{\alpha\in B_{m-l,l}}
  \Wedge^{\alpha'}(\calr^\svee)$ and we see that the latter is a
  tilting bundle by passing to the dual Grassmannian (which replaces
  $\calr^\svee$ by $\calq$) and invoking Theorem~\ref{ref-1.2-1}.
\end{proof}

For $d\ge 0$ let $\calc'_{d}$ be the category of finite-dimensional
$G_2=\GL(m-l)$-representations whose composition factors have highest
weights $\alpha$ satisfying $|\alpha|=d$ and $\alpha\in B_{m-l,l}$. We
have the following analogue of Theorem~\ref{ref-6.4-28},   where
\[
\Phi'_d\colon D^b(\calc'_d)\to \cald
\]
is defined again by $U\mapsto \call_{\GG}(U)$

\begin{theorem}
  \label{ref-6.8-34}
  The functor
$\Phi_d'$
  is fully faithful.
\end{theorem}

\begin{proof} 
  This follows by dualizing the proof of Theorem~\ref{ref-6.4-28}.
\end{proof}

Below we let $\cald'_d$ be the essential image of $D^b(\calc'_d)$
under $\Phi'_d$. We obtain the following analogue of Corollary~\ref{ref-6.5-30}.
\begin{lemma} 
  \label{ref-6.9-35}
  $\cald'_d$ is generated by $\Wedge^{\alpha'}\calr$ for $\alpha\in
  B_{m-l,l}$, $|\alpha|=d$.\qed
\end{lemma}

\begin{theorem}
  \label{ref-6.10-36}
  The triangulated category $\cald$ has a semi-orthogonal decomposition
  \[
  \left\langle\cald'_{l(m-l)},\ldots,\cald'_0\right\rangle\,.
  \]
    Furthermore $\cald'_{d}$ is
  generated by 
  $K^{\alpha}\calr$ for $\alpha\in
  B_{m-l,l}$ with $|\alpha|=d$.
\end{theorem}

\begin{proof}
  This follows by dualizing the proof of Theorem~\ref{ref-6.6-31}.
\end{proof}

\begin{remark}
  \label{rem:generators2}
  Let $\rep_e(G_i)$ be the category of finite-dimensional
  $G_i$-representations of degree $e$.  
Then  $\call_\GG(U)$ for $U \in \rep_e(G_i)$ with $e \leq d$ 
is contained in $\langle \cald_e\rangle_{e\le d}$ if $i=1$ and $\langle\cald_e'\rangle_{e\le d}$
if $i=2$. Note that
  we do not assume that the dominant weights of $U$ are in
  $B_{l,m-l}$ or $B_{m-l,l}$.
  To prove this for $i=1$ we have to show that $\Hom(\cald_f,\call_{\GG}(U))=0$
  for $f>d$. Given that $\cald_f$ is generated by
  $\Wedge^{\alpha'}\calq$ for $\alpha\in B_{l,m-l}$ and $|\alpha|=f$ this follows from
  Proposition~\ref{ref-4.3}. The argument for $i=2$ is dual.
\end{remark}

The following result finishes this section.

\begin{theorem} 
  \label{ref-6.11-37} The semi-orthogonal decompositions
  \[
  \cald=\left\langle\cald_0,\ldots,\cald_{l(m-l)}\right\rangle 
  \quad\text{and}\quad
  \cald=\left\langle\cald'_{l(m-l)},\ldots,\cald'_0\right\rangle
  \]
  are dual to each other.
  Furthermore the induced equivalence $\gamma_d\colon \cald_d\to
  \cald'_d$ defined in Lemma~\ref{ref-5.3-25} sends $L^\alpha\calq$ to
  $K^{\alpha'} \calr[d]$ for $\alpha\in B_{l,m-l}$ with $|\alpha|=d$.
\end{theorem}

\begin{proof}
  To prove that the semi-orthogonal decompositions are dual, according
  to \S\ref{ref-5-22} we have to show that
  \[
  \cald_{\le d}=\cald'_{\le d}\,,
  \]
  where we set $\cald_{\le
    d}=
  \left\langle\cald_0,\ldots,\cald_d\right\rangle$ and $\cald'_{\le
    d}=
  \left\langle\cald'_{d},\ldots,\cald'_0\right\rangle$.   We prove the
  inclusion $\cald_{\le d}\subset \cald'_{\le d}$. The opposite
  inclusion is similar.

  From Theorem~\ref{ref-6.6-31} we obtain that $\cald_{\le d}$ is
  generated by $L^\alpha\calq$ for $|\alpha|\le d$.  Thus we have to
  show that for such $\alpha$ we have $L^\alpha\calq\in \cald'_{\le
    d}$.

  According to \cite[Ch 2, Ex. 21]{Weyman:book} we have a resolution
  for $L^\alpha \calq$ given by the Schur complex
  \[
  L^\alpha(\calr\to F^\svee\otimes \calo_{\GG})
  \]
  and furthermore by \cite[Thm.
  (2.4.10)(b)]{Weyman:book} $L^\alpha(\calr\to F^\svee\otimes
  \calo_{\GG})$ has a filtration such that
  \begin{equation}
    \label{ref-6.4-38}
    \gr L^\alpha(\calr\to F^\svee\otimes \calo_{\GG})_t=
    \bigoplus_{|\nu|=t,\nu\subset\alpha}K^{\nu'}\calr\otimes
    L^{\alpha/\nu} (F^\svee)\,. 
  \end{equation}
  By Remark~\ref{rem:generators2} all $K^{\nu'}\calr$ are in $\cald'_{\le
    d}$. Hence so is $L^\alpha Q$.
  
  Assume now $|\alpha|=d$. In that case \eqref{ref-6.4-38} shows that
  \[
  L^\alpha \calq=K^{\alpha'} \calr[|\alpha|]\quad \mod \cald'_{\le
    d-1}\,.
  \]
  If in addition $\alpha\in B_{l.m-l}$ then  Lemma~\ref{ref-6.9-35} implies
  $K^{\alpha'}\calr[|\alpha|]\in \cald'_d$, from which we conclude
  that $\gamma_d(L^\alpha \calq)=K^{\alpha'}\calr[|\alpha|]$.
\end{proof}

\section{Some more comments on representation theory}
\label{sect:morecomments}
If we combine Theorems~\ref{ref-6.4-28},~\ref{ref-6.8-34},
and~\ref{ref-6.11-37} we obtain an equivalence of categories
$\delta_d\overset{\text{def}}{=}(\Phi^{\prime -1}_d\circ \gamma_d\circ
\Phi_d)[-d]$ between $D^b(\calc_d)$ and $D^b(\calc'_d)$.  The
existence of such an equivalence is well-known (see e.g.\ \cite[Cor
(3.9)]{Donkin:1993} for a similar result) but the standard
construction uses the representation theory of the symmetric group.

Below we  list some properties of the equivalence, which we will use
in \S\ref{sect:quasihereditary}.  For $\alpha$ a
partition let $M^\alpha$ be the indecomposable tilting
$G_1$-representation with highest weight $\alpha$.  Similarly for
$\beta\in B_{m-l,l}$ let $\Sigma^{\prime \beta}$ be the simple
$G_2$-representation with highest weight $\beta$ and let
$P^{\prime\beta}$ be the projective cover of $\Sigma^{\prime \beta}$
in $\calc'_d$.

Since $M^\alpha$ has highest weight $\alpha$ and since $B_{l,m-l}$ is
a saturated set of partitions, all the dominant weights of $M^\alpha$
are in $B_{l,m-l}$, whence $M^\alpha\in \calc_d$ by \cite[Lemma
E.3]{Jantzen}.
\begin{proposition}
  \label{ref-7.1-39} We have for $\alpha\in B_{l,m-l}$ with
  $|\alpha|=d$
  \begin{align}
    \delta_d(L^\alpha V)&=K^{\alpha'} W \label{ref-7.1-40}\\
    \delta_d(M^\alpha)&=P^{\prime \alpha'}\,. \label{ref-7.2-41}
  \end{align}
\end{proposition}

\begin{proof}
  Statement \eqref{ref-7.1-40} follows from Theorem~\ref{ref-6.11-37}. To prove
  \eqref{ref-7.2-41} we first note that by suitably filtering
  $M^\alpha$ and invoking \eqref{ref-7.1-40} we obtain that
  $\delta_d(M^\alpha)\in \calc'_d$. Furthermore since $\delta_d$ is an
  equivalence for $i>0$ and $\beta \in B_{l,m-l}$ with $|\beta|=d$ we have 
  \begin{multline}
    \label{ref-7.3-42}
    \Ext^i_{G_2}(\delta_d(M^\alpha),K^{\beta'}
    W)=\Ext^i_{G_2}(\delta_d(M^\alpha),\delta_d(L^{\beta} V))=
    \Ext^i_{G_1}(M^\alpha,L^\beta V)=0\,.
  \end{multline}
  We now claim that we  have for $i>0$
  \begin{equation}
    \label{ref-7.4-43}
    \Ext^i_{G_2}(\delta_d(M^\alpha),\Sigma^{\beta'})=0\,.
  \end{equation}
  We prove this by induction. As before we have a short exact sequence
  \[
  0\to U\to K^{\beta'}W\to \Sigma^{\beta'}\to 0
  \]
  where $U$ is obtained through extensions involving only
  $\Sigma^{\gamma}$ with $\gamma<\beta'$. By induction we may assume
  $\Ext^i_{G_2}(\delta_d(M^\alpha),U)=0$. Then \eqref{ref-7.4-43}
  follows from \eqref{ref-7.3-42}.  We conclude that
  $\delta_d(M^\alpha)$ is projective. Since $M^\alpha$ is
  indecomposable, the same is true for $\delta_d(M^\alpha)$.  Hence
  $\delta_d(M^\alpha)$ is equal to some $P^{\prime\gamma}$. To prove
  that $\delta_d(M^\alpha)=P^{\prime\alpha'}$ it is sufficient to
  construct a surjective map
  \[
  \delta_d(M^\alpha)\to K^{\alpha'} W
  \]
  since $K^{\alpha'}W$ has simple top $\Sigma^{\prime\alpha'}$. 
  
  By \cite[\S E.4]{Jantzen} we have a short exact sequence
  \[
  0\to H\to M^\alpha \to L^\alpha V\to 0
  \]
  where $H$ is an extension of $L^\gamma W$ with $\gamma<\alpha$.
  After applying $\delta_d$ this becomes a distinguished triangle
  \[
  \delta_d(H)\to \delta_d(M^\alpha)\to K^{\alpha'} W\to 
  \]
  with $\delta_d(H),\delta_d(M^\alpha)\in \calc'_d$. The long exact
  sequence for cohomology shows that $\delta_d(M^\alpha) \to
  K^{\alpha'} W$ is indeed surjective.
\end{proof}

\section{Exceptional sequences on Grassmannians}
\begin{proposition}
  \label{ref-8.1-44}
  Assume $\alpha,\beta\in B_{l,m-l}$ with $|\alpha|=|\beta|$. If
  \begin{equation}
    \rHom_{\calo_\GG}(L^\alpha\calq,L^\beta\calq)\neq 0
  \end{equation}
  then $\alpha\ge \beta$. Furthermore we also have
  \[
  \rHom_{\calo_\GG}(L^\alpha\calq,L^\alpha\calq)=K\,.
  \]
\end{proposition}

\begin{proof}
  We have $L^\alpha\calq=\Phi_d(L^\alpha V)$, $L^\beta\calq=
  \Phi_d(L^\beta V)$.  So to prove the first claim, by
  Theorem~\ref{ref-6.8-34} we must show that
  \[
  \rHom_{G_1}(L^\alpha V,L^\beta V)\neq 0
  \]
  implies $\alpha\ge\beta$.  Since $L^\alpha V$, $L^\beta V$ are
  induced representations it suffices to invoke \cite[Prop.\
  II.6.20]{Jantzen}.

  By \cite[Prop. II.2.8]{Jantzen} we have $\Hom_{G_1}(L^\alpha
  V,L^\alpha V)=K$. Hence to prove the second claim we have to show
  \[
  \Ext^i_{G_1}(L^\alpha V,L^\alpha V)=0
  \]
  for $i>0$. This follows from \cite[Prop.\ II.6.20]{Jantzen}.
\end{proof}

\begin{proposition}
  \label{ref-8.2-45}
  Assume $\alpha,\beta\in B_{m-l,l}$, $|\alpha|=|\beta|$. Then
  \begin{equation}
    \rHom_{\calo_\GG}(L^\alpha\calr,L^\beta\calr)\neq 0
  \end{equation}
  implies $\alpha\ge \beta$. Furthermore we also have
  \[
  \rHom_{\calo_\GG}(L^\alpha\calr,L^\alpha\calr)=K\,.
  \]
\end{proposition}

\begin{proof}
  This is proved in exactly the same way as Proposition~\ref{ref-8.1-44}.
\end{proof}

\begin{cor}
  \label{cor:exc-seqs}
  The subcollections 
  \[
  \left(L^\alpha\calq\right)_{\alpha \in B_{l,m-l},\
    |\alpha|=d}
  \quad\text{and}\quad
  \left(L^\alpha\calr\right)_{\alpha \in B_{l,m-l},\
    |\alpha|=d}
  \]
  form exceptional collections in $\cald_d$. \qed
\end{cor}

\begin{proposition} 
  \label{ref-8.3-46} For $\alpha,\beta\in B_{m-l,l}$, $|\alpha|=|\beta|$ we have
  \[
  \rHom_{\calo_\GG}(K^\alpha \calr,L^\beta \calr)=\delta_{\alpha,\beta}\cdot K\,.
  \]
\end{proposition}

\begin{proof} Put $d=|\alpha|=|\beta|$.  We have
  $K^\alpha\calr=\Phi'_d(K^\alpha W)$, $L^\beta\calr= \Phi'_d(L^\beta
  W)$.  So by Theorem~\ref{ref-6.8-34} we must show
  \[
  \rHom_{G_2}(K^\alpha W,L^\beta W)=\delta_{\alpha,\beta}\cdot K\,.
  \]
  As $K^\alpha W$ is a Weyl representation and $L^\beta W$ is an
  induced representation, it suffices to invoke
  \cite[II.4.13]{Jantzen}.
\end{proof}
Now we make $B_{l,m-l}$ into a totally ordered set by equipping it
with an arbitrary total ordering $\prec$ such that if
$|\alpha|<|\beta|$ then $\alpha\prec \beta$ and if $|\alpha|=|\beta|$
and $\alpha> \beta$ in the standard partial order on partitions, then
$\alpha\prec \beta$.  We write $\bar{B}_{l,m-l}$ for $B_{l,m-l}$,
equipped with the opposite ordering.

The following is the main result of this section.
\begin{theorem}
  \label{ref-8.4-47}
  The collections $(L^\alpha\calq)_{\alpha\in B_{l,m-l}}$ and
  $(L^{\beta'}\calr[|\beta|])_{\beta\in \bar{B}_{l,m-l}}$ form dual
  exceptional collections in $\cald$. In
  other words for $\alpha,\beta\in B_{l,m-l}$ we have 
  \begin{equation}
    \label{ref-8.3-48}
    \rHom_{\calo_\GG}(L^\alpha
    \calq,L^{\beta'}\calr[|\beta|])=\delta_{\alpha,\beta}\cdot K\,. 
  \end{equation}
\end{theorem}

\begin{proof} The fact that $(L^\alpha\calq)_{\alpha\in B_{l,m-l}}$ is
  an exceptional sequence follows from Theorem~\ref{ref-6.6-31} and
  Proposition~\ref{ref-8.1-44}. Similarly the fact that
  $(L^{\beta'}\calr[|\beta|])_{\beta\in \bar{B}_{l,m-l}}$ is an
  exceptional collection follows from Theorem~\ref{ref-6.10-36} and
  Proposition~\ref{ref-8.2-45}. So it remains to prove the duality
  property \eqref{ref-8.3-48}.  By Theorem~\ref{ref-6.11-37} combined
  with \eqref{ref-5.1-24} we may assume that $|\alpha|=|\beta|$.  We
  compute for all $i \in \ZZ$
  \begin{align*}
    \Ext^i_{\calo_\GG}(L^\alpha \calq,L^{\beta'}\calr[|\beta|])& =
    \Ext^i_{\calo_\GG}(\gamma_{|\alpha|}(L^\alpha \calq),L^{\beta'}\calr[|\beta|])\\
    &=\Ext^i_{\calo_\GG}(K^{\alpha'}
    \calr[|\alpha|],L^{\beta'}\calr[|\beta|])\\
    &=\delta_{\alpha,\beta}\cdot \delta_{i0} \cdot K\,,
  \end{align*}
  using, respectively, Lemma~\ref{ref-5.3-25}, 
  Theorem~\ref{ref-6.11-37}, and 
  Proposition~\ref{ref-8.3-46}.
\end{proof}

To conclude this section we use the ``linkage principle''
\cite[Cor. II.6.17]{Jantzen} to recover the result of Kaneda, mentioned
in the Introduction, that Kapranov's tilting bundle $\calk$ remains
tilting in large characteristic. 

\begin{lemma} \label{lem:linkage} 
  Assume that $K$ has characteristic
  $p$ with $p \geq m-1$.  Let $\alpha\in B_{l,m-l}$. Then
  $\L^{\alpha'}\calq$ is a direct sum of $L^\beta \calq$ with
  $|\beta|=|\alpha|$ and furthermore there are no homomorphisms between
  the summands of $\L^{\alpha'}\calq$.
\end{lemma}

\begin{proof} Set $d=|\alpha|$.  Using Theorem~\ref{ref-6.4-28} it is
  enough to prove the following claim: $\calc_d$ is a semi-simple
  category with simple objects given by $L^\beta V$ for $\beta\in
  B_{l,m-l}$ with $|\beta|=d$.  Indeed if this claim holds then
  $\L^{\alpha'} V$ is a direct sum of the simple objects $L^\beta V$
  (which thus have no $\Hom$'s among them) and it suffices to apply
  the fully faithful functor $\Phi(-)_d$ to obtain the corresponding
  result for $\L^{\alpha'}\calq$.

  The claim follows directly from the linkage principle
  \cite[Cor. II.6.17]{Jantzen} which we state in the case of interest
  to us.  If $\gamma$, $\delta$ are dominant weights for $G_1$ and
  $\Ext^1_{\calo_{\GG}}(\Sigma^\gamma, \Sigma^\delta)\neq 0$ then
  $\gamma$, $\delta$ are in the same orbit for the \emph{affine Weyl
    group}.


  A fundamental domain\footnote{%
    A fundamental domain is a complete irredundant set of
    orbit representatives~\cite[IV.3.3]{Bourbaki:Lie4-6}.} 
  $\overline C$ for the affine Weyl group
  (\cite[II.6.1(6)]{Jantzen}) is given by 
  the set of $x= [x_1, \dots, x_l]$ satisfying $\langle x+\rho,
  \alpha\rangle \leq p$ for all positive roots $\alpha$, where
  $\rho=[l, l-1, \dots, 1]$; equivalently
  \[
  0 \leq x_i-i-x_j+j \leq p
  \]
  for $j>i$.  The first inequality is automatically satisfied for a
  dominant weight. For the second inequality we note that if
  $\gamma = [\gamma_1, \dots, \gamma_l]\in B_{l,m-l}$ then
  \[
  \gamma_i-\gamma_j\le m-l
  \]
  and 
  \[
  -i+j\le l-1\,.
  \]
  Thus
  \[
  \gamma_i-i-\gamma_j+j \leq m-l+l-1= m-1\le p\,.
  \]
  In other words
  \(
  B_{l,m-l}\subseteq \overline C
  \)
  and thus no two elements of $B_{l,m-l}$ are in the same orbit for
  the affine Weyl group. The claim follows.
\end{proof}

The fact that $\bigoplus_{\alpha\in B_{l,m-l}} \L^{\alpha'}\calq$ is a
tilting object, together with Theorem~\ref{ref-6.6-31} and the
previous lemma, yields immediately the following.

\begin{corollary}
   [(Kaneda~\cite{Kaneda:2008})] 
  \label{cor:kaneda}
  The Kapranov strong exceptional collection
  $(L^\alpha\calq)_{\alpha\in B_{l,m-l}}$ remains strong exceptional
  as long as $p\ge m-1$.  In particular $\calk = \bigoplus_{\alpha\in
    B_{l,m-l}} L^\alpha\calq$ remains tilting for such $p$.
\end{corollary}

\section{Relation with quasi-hereditary algebras}
\label{sect:quasihereditary}

We quickly remind the reader of the module-theoretic description of
quasi-hereditary algebras \`a la Dlab-Ringel~\cite{Dlab-Ringel:1992};
see also for example~\cite{Erdmann:1994, Hille-Perling:2011arxiv}.
Let $A$ be a finite-dimensional $K$-algebra and let
$(S(\lambda))_{\lambda \in \Lambda}$ be a complete set of the simples,
with projective covers $P(\lambda) \onto S(\lambda)$ and injective
hulls $S(\lambda) \into Q(\lambda)$.

Fix a total ordering $\prec$ on $\Lambda$. Define the \emph{standard
  module} $\Delta(\lambda)$ to be the largest quotient of $P(\lambda)$
having composition factors of the form $S(\mu)$ with $\mu \preceq
\lambda$.  Equivalently~\cite[Lemma 1.1]{Dlab-Ringel:1992}
$\Delta(\lambda)$ is the quotient of $P(\lambda)$ by the maximal
submodule generated by any direct sum of the form $\bigoplus_{\mu
  \succ \lambda} P(\mu)$.  Similarly the \emph{costandard module}
$\nabla(\lambda)$ is the largest submodule of $Q(\lambda)$ having
composition factors $S(\mu)$ with $\mu \preceq \lambda$.  In
particular $S(\lambda)$ is the top of $\Delta(\lambda)$ and the socle
of $\nabla(\lambda)$.  
Set $\Delta = (\Delta(\lambda))_{\lambda \in \Lambda}$ and $\nabla =
(\nabla(\lambda))_{\lambda \in \Lambda}$. 

We assume that each $\Delta(\lambda)$ (equivalently
each $\nabla(\lambda)$) is Schurian, i.e.\ the endomorphism ring is a
division ring.

For an arbitrary collection $\Theta$ of $A$-modules, denote by
$\calf(\Theta)$ the class of $\Theta$-filtered modules,
that is, $A$-modules $M$ having a filtration $M = M_0 \supset M_1
\supset \cdots \supset M_t = 0$ with successive quotients
$M_i/M_{i-1}$ in the collection $\Theta$.

\begin{definition}
  [({\cite[Theorem 1]{Dlab-Ringel:1992}}; see also~\cite{Donkin:1981,
    Scott:1987, Cline-Parshall-Scott:1988})] 
  The algebra $A$ (with the fixed order $\prec$) is called
  \emph{quasi-hereditary} if the following equivalent conditions hold:
  \begin{enumerate}[\quad(i)]
  \item ${}_A A \in \calf(\Delta)$;
  \item $\calf (\Delta) = \left\{X \ \middle|\
      \Ext_A^1(X,\nabla(\lambda))=0\right\}$; 
  \item $\calf (\Delta) = \left\{X \ \middle|\
      \Ext_A^i(X,\nabla(\lambda))=0 \text{ for all } i \geq 1\right\}$;
  \item $\Ext_A^2(\Delta, \nabla)=0$.
  \end{enumerate}
\end{definition}

For a quasi-hereditary algebra, the standard and costandard modules
determine each other in $D^b_f(A)$ (the bounded derived category with
finite cohomology) by
\begin{equation}
  \label{ref-9.5-54}
  \rHom_{A}(\Delta(\lambda),\nabla(\mu))=\delta_{\lambda,\mu}\cdot K\,.
\end{equation}

For the proof of Theorem~\ref{ref-9.1-49} below we also remind the
reader of the notion of
``standardization''~\cite[\S3]{Dlab-Ringel:1992}, see also the
``universal extensions'' of~\cite{Hille-Perling:2011arxiv}, in the
special case of interest.  An indexed collection $\Theta =
(\Theta(\lambda))_{\lambda \in \Lambda}$ of objects in a $K$-linear
abelian category $\calc$ is \emph{standardizable} provided
\begin{enumerate}[\quad(i)]
\item $\Hom_\calc(\Theta(\lambda),\; \Theta(\mu))$ and
  $\Ext_\calc^1(\Theta(\lambda)\; \Theta(\mu))$ are finite-dimensional
  for all $\lambda$, $\mu$, and
\item the quiver with vertex set $\Lambda$ and an arrow $\lambda \to
  \mu$ if either there is a non-trivial non-isomorphism
  $\Theta(\lambda) \to \Theta(\mu)$ or
  $\Ext_\calc^1(\Theta(\lambda),\; \Theta(\mu)) \neq 0$ has no
  oriented cycles.
\end{enumerate}
In particular note that if $\Theta$ is an exceptional
collection then $\Theta$ is standardizable.

The next result is essentially contained in the proof of~\cite[Theorem
2]{Dlab-Ringel:1992} for modules over a finite-dimensional algebra;
see also~\cite[Theorem 5.1]{Hille-Perling:2011arxiv} for a statement
in the geometric context.

\begin{theorem}
  \label{thm:DRthm2}
  Let $\Theta = (\Theta(\lambda))_{\lambda\in \Lambda}$ be a
  standardizable collection in an abelian category $\calc$. Then there exists a
  projective generator $P \in \calf(\Theta)$ such that $A'
  = \End_\calc(P)$ is quasi-hereditary with standard modules
  $\Hom_\calc(P,\Theta(\lambda))$.
\end{theorem}

As in the introduction put
  \[
  \calt=\bigoplus_{\alpha\in B_{l,m-l}} \Wedge^{\alpha'}\calq
  \]
  and $A=\End_{\calo_\GG}(\calt)$.  Denote by $A^\circ$ the opposite
  algebra.  For $\alpha\in B_{l,m-l}$ consider the following complexes
  of right $A$-modules:
  \begin{align}
    \Delta(\alpha)&=\rHom_{\calo_\GG}(\calt,L^\alpha\calq)\label{ref-9.1-50}\\
    \nabla(\alpha)&=\rHom_{\calo_\GG}(\calt,L^{\alpha'}\calr[|\alpha|])\label{ref-9.2-51}\\
    P(\alpha)&=\rHom_{\calo_\GG}(\calt,\call_\GG(M^\alpha))\label{ref-9.3-52}\\
    S(\alpha)&=\rHom_{\calo_\GG}(\calt,\call_\GG(\Sigma^{\prime\alpha'})[|\alpha|])\label{ref-9.4-53}
  \end{align}
  where as before $M^\alpha$ is the indecomposable tilting
   $\GL(l)$-representation with highest weight $\alpha$.

\begin{theorem}
  \label{ref-9.1-49}\ 
  The complexes $\Delta(\alpha)$, $\nabla(\alpha)$, $P(\alpha)$,
    $S(\alpha)$ are concentrated in degree zero, and the
    $S(\alpha)$ are the simple right $A$-modules with the $P(\alpha)$
    their projective covers.
  
  As in \S7, let $\prec$ be an arbitrary total ordering on
    $B_{l,m-l}$ such that if $|\alpha| < |\beta|$, or if
    $|\alpha|=|\beta|$ and $\alpha > \beta$ in the natural partial
    order on partitions, then $\alpha \prec \beta$.  Then $A$ is
    quasi-hereditary with respect to this ordering. The standard and
    costandard modules having $S(\alpha)$ respectively as top and
    socle are $\Delta(\alpha)$ and $\nabla(\alpha)$.
\end{theorem}

\begin{proof}
  Set $\Theta=(L^\alpha \calq)_{\alpha\in B_{l,m-l}}$.  Then $\Theta$
  is an exceptional collection 
  by Theorem~\ref{ref-8.4-47}, so in particular is standardizable.
  Let $P \in \calf(\Theta)$ be the projective generator guaranteed by
  Theorem~\ref{thm:DRthm2}, so that $A' = \End_{\calo_\GG}(P)$ is
  quasi-hereditary with standard modules $\Hom_{\calo_\GG}(P, L^\alpha
  \calq)$.

  On the other hand Proposition~\ref{ref-1.3-3} implies that $\calt$
  is a projective generator for $\calf(\Theta)$ as well.  This easily
  yields that $A$ and $A'$ are Morita equivalent and that the objects
  $\Hom_{\calo_\GG}(\calt,L^\alpha \calq)$ are the standard objects.
  By Proposition~\ref{ref-1.3-3} we may replace $\Hom$ by $\rHom$.

  Since the costandard modules $\nabla(\beta)$ are characterized by
  $\rHom_{A^\circ}(\Delta(\alpha),\;
  \nabla(\beta))=\delta_{\alpha,\beta}\cdot K$, we deduce from
  \eqref{ref-8.3-48} that they are indeed given by the formula
  \eqref{ref-9.2-51}.

  By~\cite[Lemma (3.4)]{Donkin:1993} the indecomposable summands of
  $\Wedge^{\alpha'}V$ are precisely the tilting representations
  $M^\beta$, and $M^\alpha$ occurs as a direct summand with
  multiplicity one in $\Wedge^{\alpha'} V$. It follows that the
  $P(\alpha)$ are indeed the indecomposable projectives.

  To show that the $S(\alpha)$ are the corresponding simple
  $A$-modules we have to prove
  \[
  \rHom_{A^\circ}(P(\alpha),S(\beta))=\delta_{\alpha,\beta}\cdot K\,.
  \]
  We compute
  \begin{align*} \rHom_{A^\circ}(P(\alpha), S(\beta)) & =
    \rHom_{\calo_\GG}(\call_\GG(M^\alpha), \call_\GG(\Sigma^{\prime\beta'})[|\beta|])\,. 
  \end{align*}
  It follows from Theorem~\ref{ref-6.10-36} combined with
  \eqref{ref-5.1-24} that if $|\alpha|\neq |\beta|$ then there is
  nothing to prove. So assume $|\alpha|=|\beta|=d$.  Then we have
  \begin{align*}
    \rHom_{\calo_\GG}(\call_\GG(M^\alpha), \call_\GG(\Sigma^{\prime\beta'})[|\beta|])
    &=\rHom_{\calo_\GG}(\gamma_{d}(\Phi_d(M^\alpha)), \Phi'_d(\Sigma^{\prime\beta'})[|\beta|])\\ 
    &=\rHom_{G_2}(\delta_d(M^\alpha),\Sigma^{\prime\beta'})\\
    &=\rHom_{G_2}(P^{\prime\alpha'},\Sigma^{\prime\beta'})\\
    &=\delta_{\alpha,\beta}\cdot K\,,
  \end{align*}
  where the first equality is Lemma~\ref{ref-5.3-25} and the third is
  Proposition~\ref{ref-7.1-39}; in the second, $\delta_d$ is as introduced
  in \S\ref{sect:morecomments}. 

  It remains to prove that $S(\alpha)$ is the top of $\Delta(\alpha)$
  and the socle of $\nabla(\alpha)$. Since for quasi-hereditary
  algebras the top of $\Delta(\alpha)$ coincides with the socle of
  $\nabla(\alpha)$ it is sufficient to prove only the first of these
  statements.  There is a surjective map $M^\alpha\to L^\alpha V$
  whose kernel is an extension of $L^\beta V$ with $\beta<\alpha$
  \cite[\S E.4]{Jantzen}.  Apply
  $\Hom_{\calo_\GG}(\calt,\call_\GG(-))$ to obtain a map $P(\alpha)
  \to \Delta(\alpha)$, which is surjective by
  Proposition~\ref{ref-1.3-3}. This finishes the proof.
\end{proof}

\begin{example}
  \label{eg:G24}
  We compute the quiver and relations of the quasi-hereditary algebra
  $A$ in the first non-trivial example, $(m,l)=(4,2)$.  We live
  inside the $2 \times 2$ box $B_{2,2}$, so the vertices of the
  quiver, equivalently the summands of the tilting bundle $\calt$, are
  labeled 
  \[
  \calo,\ \calq,\ \Wedge^2\calq,\ \calq\otimes\calq,\ 
  \Wedge^2\calq\otimes\calq,\ \left(\Wedge^2\calq\right)^{\otimes 2}\,.
  \]
  The quiver has the following form.
  \[
  \xymatrix{
    & & \Wedge^2 \calq \ar@<0.2pc>[dd]^a\\
    \calo \ar[r]^{s_\lambda}  & \calq  \ar[dr]_{\alpha_\lambda}
        & & \Wedge^2\calq\otimes\calq \ar[r]^{t_\lambda}&
       \left(\Wedge^2\calq\right)^{\otimes 2} \\
     & & \calq\otimes\calq \ar@<0.2pc>[uu]^p  \ar@<-0.2pc>[ur]_{\beta_\lambda}
  }
  \]
  The labels stand for
  natural maps between these bundles, some of which depend on a global
  section $\lambda \in F^\svee$:
  \begin{enumerate}[\quad$\bullet$]
  \item $p \colon \calq\otimes \calq \to \Wedge^2 \calq$ the natural
    surjection and $a \colon \Wedge^2 \calq \to \calq \otimes \calq$ the anti-symmetrization;
  \item $s_\lambda\colon \calo \to \calq$ with $1 \mapsto \lambda$;
  \item $\alpha_\lambda \colon \calq \to \calq\otimes \calq$ with $x
    \mapsto \lambda \otimes x$\,;
  \item $\beta_\lambda \colon \calq \otimes \calq \to \Wedge^2\calq
    \otimes \calq$ with $x\otimes y \mapsto \lambda \wedge y \otimes x$\,; and
  \item $t_\lambda \colon \Wedge^2\calq \otimes \calq \to \left(\Wedge^2
    \calq\right)^{\otimes 2}$ with $x\wedge y \otimes z \mapsto
    x\wedge y \otimes \lambda \wedge z$.
  \end{enumerate}
  These maps generate all the arrows in the quiver.  For example, the
  obvious complementary map $\alpha'_\lambda \colon \calq \to \calq
  \otimes \calq$ defined by $\alpha'_\lambda (x) = x \otimes \lambda$
  can be obtained as $\alpha_\lambda(1-ap)$.  The relations are most
  compactly written in terms of the pseudo-idempotent $e
  \overset{\text{def}}= ap$ satisfying $e^2=2e$, and the ``swap'' $1-e$ which sends $x
  \otimes y$ to $y \otimes x$. We have
  \begin{enumerate}[\quad$\bullet$] 
  \item $pa = 2\Id_{\bigwedge^2 \calq}$\,;
  \item $(1-e)\alpha_\lambda s_\mu = \alpha_\mu s_\lambda$\,;
  \item $\beta_\lambda(1-e)\alpha_\mu = \beta_\lambda \alpha_\mu -
    \beta_\mu \alpha_\lambda$\,;
  \item $t_\lambda \beta_\mu (1-e) = t_\mu \beta_\lambda$\,.
  \end{enumerate}

  We observe that in this picture, each vertical ``slice'' is
  equivalent to the derived category of a generalized Schur algebra.
  For example, in the middle we recognize the quiver for the Schur
  algebra $S(2,2)$ in characteristic $2$~\cite[3.1.1,
  5.4]{Erdmann:1993}. 

  In characteristic different from $2$, the idempotent $\frac 12 e$
  gives $\calq \otimes \calq
  \cong \Wedge^2 \calq \oplus \Sym^2 \calq$, and the algebra
  becomes Morita-equivalent to the path algebra of the
  equivariant quiver
  \[
  \xymatrix{
    & & \Wedge^2 \calq \ar[dr]^{F^\svee}  \ar@/^1.5pc/@{-->}[drr]^{D^2
      F^\svee}\\
    \calo \ar[r]^{F^\svee} \ar@/^1.5pc/@{-->}[urr]^{D^2F^\svee} 
    \ar@/_1.5pc/@{-->}[drr]_{\bigwedge^2F^\svee}
      & \calq \ar[dr]_{F^\svee} \ar[ur]^{F^\svee}
      \ar@{-->}[rr]^{F^\svee \otimes F^\svee}
       & & \Wedge^2\calq\otimes\calq \ar[r]^{F^\svee} 
       &
       \left(\Wedge^2\calq\right)^{\otimes2} \\
    & & \Sym^2\calq  \ar[ur]_{F^\svee} 
       \ar@/_1.5pc/@{-->}[urr]_{\bigwedge^2F^\svee}
  }
  \]   
  with relations
  \begin{enumerate}[\quad$\bullet$]
  \item $\calo \to \Wedge^2\calq$ given by $D^2 F^\svee$;
  \item $\calo \to \Sym^2 \calq$ given by $\Wedge^2 F^\svee$;
  \item $\calq \to \Wedge^2 \calq \otimes \calq$ given by $F^\svee \otimes F^\svee$;
  \item $\Wedge^2 \calq \to \left(\Wedge^2 \calq\right)^{\otimes 2}$
    given by $D^2 F^\svee$; and
  \item $\Sym^2 \calq \to \left(\Wedge^2 \calq\right)^{\otimes 2}$
    given by $\Wedge^2 F^\svee$.
  \end{enumerate}
  Most of these are straightforward to verify.  The relations across
  the central diamond, however, are not the obvious commutativity
  ones~\cite{Hille:1998, Sam-Weyman:2011}. To compute those relations, give names to
  the maps:
  \[
  \xymatrix{
    & \Wedge^2 \calq \ar[dr]^{b_\lambda}  \\
      \calq \ar[dr]_{c_\lambda} \ar[ur]^{a_\lambda}
       & & \Wedge^2\calq\otimes\calq  \\
    & \Sym^2\calq  \ar[ur]_{d_\lambda} 
  }
  \]  
  with
  \begin{enumerate}[\quad$\bullet$]
  \item $a_\lambda(x) = \lambda \wedge x$;
  \item $b_\lambda(x\wedge y) = x\wedge y \otimes \lambda$;
  \item $c_\lambda(x) = \lambda x$; and
  \item $d_\lambda(xy) = \lambda \wedge x \otimes y + \lambda \wedge y
    \otimes x$.
  \end{enumerate}
  Then we find, for $\lambda, \mu \in F^\svee$,
  \[
  d_\mu c_\lambda =
  2 b_\lambda a_\mu - b_\mu a_\lambda\,.
  \]
  It follows that the defining relations are
  \begin{align*}
     d_\lambda c_\mu + d_\mu c_\lambda 
       &= b_\lambda a_\mu + b_\mu a_\lambda \\
    d_\mu c_\lambda - d_\lambda c_\mu
       &= 3(b_\lambda a_\mu - b_\mu a_\lambda)\,.
  \end{align*}
\end{example}


\begin{thebibliography}{BLVdB13}

\bibitem[BK89]{Bondal-Kapranov:1989}
Alexei~I. Bondal and Mikhail~M. Kapranov, \emph{Representable functors, {S}erre
  functors, and reconstructions}, Izv. Akad. Nauk SSSR Ser. Mat. \textbf{53}
  (1989), no.~6, 1183--1205, 1337. \MR{1039961}

\bibitem[BLVdB13]{BLVdB2}
Ragnar-Olaf Buchweitz, Graham~J. Leuschke, and Michel Van~den Bergh,
  \emph{Non-commutative desingularization of determinantal varieties, {II}:
  {A}rbitrary minors}, preprint (2013), 61 pages, \url{http://arxiv.org/abs/1106.1833}.

\bibitem[Bof88]{Boffi:1988}
Giandomenico Boffi, \emph{The universal form of the {L}ittlewood-{R}ichardson
  rule}, Adv. in Math. \textbf{68} (1988), no.~1, 40--63. \MR{931171}

\bibitem[Bon89]{Bondal:1989}
Alexei~I. Bondal, \emph{Representations of associative algebras and coherent
  sheaves}, Izv. Akad. Nauk SSSR Ser. Mat. \textbf{53} (1989), no.~1, 25--44.
  \MR{992977}

\bibitem[Bou02]{Bourbaki:Lie4-6}
Nicolas Bourbaki, \emph{Lie groups and {L}ie algebras. {C}hapters 4--6},
  Elements of Mathematics (Berlin), Springer-Verlag, Berlin, 2002, Translated
  from the 1968 French original by Andrew Pressley. \MR{1890629}

\bibitem[BVdB03]{BVdB}
Alexei~I. Bondal and Michel Van~den Bergh, \emph{Generators and
  representability of functors in commutative and noncommutative geometry},
  Mosc. Math. J. \textbf{3} (2003), no.~1, 1--36, 258. \MR{1996800}

\bibitem[CPS83]{Cline-Parshall-Scott:1983}
Edward T.~Cline, Brian Parshall, and Leonard Scott, \emph{A {M}ackey imprimitivity
  theory for algebraic groups}, Math. Z. \textbf{182} (1983), no.~4, 447--471.
  \MR{701363}

\bibitem[CPS88]{Cline-Parshall-Scott:1988}
\bysame,
  \emph{Finite-dimensional algebras and highest weight categories}, J. Reine
  Angew. Math. \textbf{391} (1988), 85--99. \MR{961165}

\bibitem[Dan96]{Danilov:1996}
Vladimir~Ivanovich Danilov, \emph{Cohomology of algebraic varieties}, Algebraic
  geometry, {II}, Encyclopaedia Math. Sci., vol.~35, Springer, Berlin, 1996,
  pp.~1--125, 255--262. \MR{1392958}

\bibitem[Don81]{Donkin:1981}
Stephen Donkin, \emph{A filtration for rational modules}, Math. Z. \textbf{177}
  (1981), no.~1, 1--8. \MR{611465}

\bibitem[Don93]{Donkin:1993}
\bysame, \emph{On tilting modules for algebraic groups}, Math. Z.
  \textbf{212} (1993), no.~1, 39--60. \MR{1200163}

\bibitem[DR92]{Dlab-Ringel:1992}
Vlastimil Dlab and Claus~Michael Ringel, \emph{The module theoretical approach
  to quasi-hereditary algebras}, Representations of algebras and related topics
  ({K}yoto, 1990), London Math. Soc. Lecture Note Ser., vol. 168, Cambridge
  Univ. Press, Cambridge, 1992, pp.~200--224. \MR{1211481}

\bibitem[DRS74]{Doubilet-Rota-Stein:1974}
Peter Doubilet, Gian-Carlo Rota, and Joel Stein, \emph{On the foundations of
  combinatorial theory. {IX}. {C}ombinatorial methods in invariant theory},
  Studies in Appl. Math. \textbf{53} (1974), 185--216. \MR{0498650}

\bibitem[Erd93]{Erdmann:1993}
Karin Erdmann, \emph{Schur algebras of finite type}, Quart. J. Math. Oxford
  Ser. (2) \textbf{44} (1993), no.~173, 17--41. \MR{1206201}

\bibitem[Erd94]{Erdmann:1994}
\bysame, \emph{Symmetric groups and quasi-hereditary algebras},
  Finite-dimensional algebras and related topics ({O}ttawa, {ON}, 1992), NATO
  Adv. Sci. Inst. Ser. C Math. Phys. Sci., vol. 424, Kluwer Acad. Publ.,
  Dordrecht, 1994, pp.~123--161. \MR{1308984}

\bibitem[Ful97]{Fulton:1997}
William Fulton, \emph{Young tableaux}, London Mathematical Society Student
  Texts, vol.~35, Cambridge University Press, Cambridge, 1997, With
  applications to representation theory and geometry. \MR{1464693}

\bibitem[Hil98]{Hille:1998}
Lutz Hille, \emph{Homogeneous vector bundles and {K}oszul algebras}, Math.
  Nachr. \textbf{191} (1998), 189--195. \MR{1621314}

\bibitem[HP11]{Hille-Perling:2011arxiv}
Lutz Hille and Markus Perling, \emph{{Tilting Bundles on Rational Surfaces and
  Quasi-Hereditary Algebras}}, preprint (2011), 15 pages, \url{http://arxiv.org/abs/1110.5843}.

\bibitem[Jan03]{Jantzen}
Jens~Carsten Jantzen, \emph{Representations of algebraic groups}, second ed.,
  Mathematical Surveys and Monographs, vol. 107, American Mathematical Society,
  Providence, RI, 2003. \MR{2015057}

\bibitem[Kan08]{Kaneda:2008}
Masaharu Kaneda, \emph{Kapranov's tilting sheaf on the {G}rassmannian in
  positive characteristic}, Algebr. Represent. Theory \textbf{11} (2008),
  no.~4, 347--354. \MR{2417509}

\bibitem[Kap88]{Kapranov:1988}
Mikhail~M. Kapranov, \emph{On the derived categories of coherent sheaves on
  some homogeneous spaces}, Invent. Math. \textbf{92} (1988), no.~3, 479--508.
  \MR{939472}

\bibitem[Kem76]{Kempf:1976}
George~R. Kempf, \emph{Linear systems on homogeneous spaces}, Ann. of Math. (2)
  \textbf{103} (1976), no.~3, 557--591. \MR{0409474}

\bibitem[{Kuz}09]{Kuznetsov:2009}
Alexander {Kuznetsov}, \emph{{Hochschild homology and semiorthogonal
  decompositions}}, preprint (2009), 30 pages, \url{http://arxiv.org/abs/0904.4330}.

\bibitem[LSW89]{Levine-Srinivas-Weyman:1989}
Marc Levine, V.~Srinivas, and Jerzy Weyman, \emph{{$K$}-theory of twisted
  {G}rassmannians}, $K$-Theory \textbf{3} (1989), no.~2, 99--121. \MR{1029954}

\bibitem[Sco87]{Scott:1987}
Leonard~L. Scott, \emph{Simulating algebraic geometry with algebra. {I}. {T}he
  algebraic theory of derived categories}, The {A}rcata {C}onference on
  {R}epresentations of {F}inite {G}roups ({A}rcata, {C}alif., 1986), Proc.
  Sympos. Pure Math., vol.~47, Amer. Math. Soc., Providence, RI, 1987,
  pp.~271--281. \MR{933417}

\bibitem[SW11]{Sam-Weyman:2011}
Steven V Sam and Jerzy Weyman, \emph{Pieri resolutions for classical
  groups}, J.\ Algebra \textbf{329} (2011), 222--259. \MR{2769324}
See also the revised and augmented version at
\url{http://arxiv.org/abs/0907.4505v5}. 

\bibitem[Wey03]{Weyman:book}
Jerzy Weyman, \emph{Cohomology of vector bundles and syzygies}, Cambridge
  Tracts in Mathematics, vol. 149, Cambridge University Press, Cambridge, 2003.
  \MR{1988690}

\end{thebibliography}

\newcommand{\arxiv}[2][AC]{\mbox{\href{http://arxiv.org/abs/#2}{\textsf{arXiv:#2}}}}
\newcommand{\oldarxiv}[2][AC]{\mbox{\href{http://arxiv.org/abs/math/#2}{\textsf{arXiv:math/#2[math.#1]}}}}
\providecommand{\MR}[1]{\mbox{\href{http://www.ams.org/mathscinet-getitem?mr=#1}{#1}}}
\renewcommand{\MR}[1]{%
  {\href{http://www.ams.org/mathscinet-getitem?mr=#1}{MR #1}.}}
\providecommand{\bysame}{\leavevmode\hbox to3em{\hrulefill}\thinspace}
\newcommand{\arXiv}[1]{%
  \relax\ifhmode\unskip\space\fi\href{http://arxiv.org/abs/#1}{arXiv:#1}}

\def\cprime{$'$}
\providecommand{\bysame}{\leavevmode\hbox to3em{\hrulefill}\thinspace}
\providecommand{\MR}{\relax\ifhmode\unskip\space\fi MR }
\providecommand{\MRhref}[2]{%
  \href{http://www.ams.org/mathscinet-getitem?mr=#1}{#2}
}
\providecommand{\href}[2]{#2}

\end{document}